\documentclass{amsart}

\usepackage{amssymb,amsmath,amscd,amsthm}
\usepackage{amsfonts,epsfig,latexsym,graphicx,amssymb, pict2e}

\usepackage[T1]{fontenc}
\usepackage{color}

\date{\today}

\newcommand{\Z}{{\mathbb Z}}
\newcommand{\R}{{\mathbb R}}
\newcommand{\C}{{\mathbb C}}
\newcommand{\N}{{\mathbb N}}
\newcommand{\T}{{\mathbb T}}

\renewcommand{\S}{{\mathbb S}}

\def\a{{\sf a}}
\def\b{{\sf b}}
\def\ab{{\a\b}}
\def\alphabet{{\sf A}}
\newcommand{\dd}{{\mathrm{d}}}

\newcommand{\Leb}{{\mathrm{Leb}}}

\newcommand{\Per}{{\mathrm{Per}}}
\newcommand{\loc}{{\mathrm{loc}}}
\newcommand{\cs}{{\mathrm{cs}}}
\newcommand{\cu}{{\mathrm{cu}}}
\newcommand{\sss}{{\mathrm{ss}}}
\newcommand{\uu}{{\mathrm{uu}}}
\newcommand{\const}{{\mathrm{const}}}

\newtheorem{theorem}{Theorem}[section]
\newtheorem{lemma}[theorem]{Lemma}
\newtheorem{prop}[theorem]{Proposition}
\newtheorem{definition}[theorem]{Definition}

\theoremstyle{definition}
\newtheorem{remark}[theorem]{Remark}

\theoremstyle{definition}
\newtheorem{defi}[theorem]{Definition}

\theoremstyle{definition}

\theoremstyle{definition}

\theoremstyle{definition}
\newtheorem{quest}{Question}

\theoremstyle{definition}
\newtheorem{assumption}[theorem]{Assumption}

\allowdisplaybreaks

\numberwithin{equation}{section}

\newcommand{\tr}{\mathrm{tr} }

\newcommand{\eqdef}{\overset{\mathrm{def}}=}

\begin{document}

\title[Spectra of Multidimensional Schr\"odinger Operators]{Multidimensional Schr\"odinger Operators Whose Spectrum Features a Half-Line and a Cantor Set}

\author[D.\ Damanik]{David Damanik}
\address{Department of Mathematics, Rice University, Houston, TX~77005, USA}
\email{damanik@rice.edu}
\thanks{D.D.\ was supported in part by NSF grant DMS--1700131 and by an Alexander von Humboldt Foundation research award}

\author[J.\ Fillman]{Jake Fillman}
\address{Department of Mathematics, Texas State University, San Marcos, TX~78666}
\email{fillman@txstate.edu}

\author[A.\ Gorodetski]{Anton Gorodetski}
\address{Department of Mathematics, University of California, Irvine, CA~92697, USA\\
and  National Research University Higher School of Economics, Russian Federation }
\email{asgor@uci.edu}
\thanks{A.G.\ was supported in part by Simons Fellowship (grant number 556910), NSF grant DMS--1855541, and by Laboratory of Dynamical Systems and Applications NRU HSE, grant of the Ministry of science and higher education of the RF ag. N 075-15-2019-1931.}


\begin{abstract}
We construct multidimensional Schr\"odinger operators with a spectrum that has no gaps at high energies and that is nowhere dense at low energies. This gives the first example for which this widely expected topological structure of the spectrum in the class of uniformly recurrent Schr\"odinger operators, namely the coexistence of a half-line and a Cantor-type structure, can be confirmed. Our construction uses Schr\"odinger operators with separable potentials that decompose into one-dimensional potentials generated by the Fibonacci sequence and relies on the study of such operators via the trace map and the Fricke-Vogt invariant. To show that the spectrum  contains a half-line, we prove an abstract Bethe--Sommerfeld criterion for sums of Cantor sets which may be of independent interest.
\end{abstract}

\maketitle



\section{Introduction}

The spectral analysis of Schr\"odinger operators
$$
H_V = - \Delta + V \quad \text{in } L^2(\R^d)
$$
plays a central role in quantum mechanics. The most fundamental issue is the study of the spectrum of $H_V$,
$$
\sigma(H_V) = \{ E : (H_V-E)^{-1} \text{ does not exist as a bounded operator} \}.
$$

The simplest case, $V \equiv 0$, can be readily understood via the Fourier transform and one finds $\sigma(H_0) = [0,\infty)$. The classical Weyl theorem then shows that $\sigma_\mathrm{ess}(H_V) = [0,\infty)$ as long as $V$ is a relatively compact perturbation of $H_0 = -\Delta$. A sufficient condition for that is, e.g., $V \in L^p(\R^d) + L^\infty_\varepsilon(\R^d)$ (i.e.\ for every $\varepsilon > 0$, we can write $V = V_{1,\varepsilon} + V_{2,\varepsilon}$ with $V_{1,\varepsilon} \in L^p(\R^d)$ and $\|V_{2,\varepsilon}\|_\infty < \varepsilon$) with $p \ge \max \{2, d/2\}$ if $d \not= 4$ and $p > 2$ if $d = 4$; compare \cite[Problem XIII.41]{RS4}. In particular, $\sigma(H_V)$ contains a half line if $V$ is (locally sufficiently regular and) arbitrarily small near infinity.

On the other hand, the spectrum may have gaps at high energies when $V$ is not arbitrarily small near infinity. For an explicit example, one may consider the case $d = 1$ and $V(x) = \cos x$; see \cite[Example~1 in Section XIII.16]{RS4}. More generally, a generic periodic potential $V$ will give rise in one space dimension to a Schr\"odinger operator whose spectrum does not contain a half line \cite{S1976}. However, this is a strictly one-dimensional phenomenon. If $d \ge 2$ and $V$ is periodic (i.e., the periodicity vectors of $V$ form a lattice), then $\sigma(H_V)$ does contain a half line. This statement is known as the \textit{Bethe-Sommerfeld conjecture}, which is in fact a theorem due to work of several authors; see \cite{Parn2008AHP, Skr79, Skr84, Skr85, Vel88} and references therein. In recent years, there has also been interest in analogs of the Bethe--Sommerfeld conjecture for discrete Schr\"odinger operators on periodic lattices \cite{EF2017,FH,HanJit,KrugPreprint} and quantum graphs~\cite{ET2017}.

It is natural to take this further and attempt to establish the Bethe-Sommerfeld property (i.e., the presence of a half line in the spectrum) for more general potentials $V$. On the positive side, there has been recent work establishing this result for certain classes of almost periodic potentials in dimension two \cite{KL13, KS18+}. However, there is also a negative result: it was shown in \cite{DFG2019} that there are almost periodic potentials $V$ in any dimension for which $\sigma(H_V)$ is a generalized Cantor set, that is, a perfect set with empty interior. This result provides the first example of an almost periodic Schr\"odinger operator in dimension $d \ge 2$ for which the spectrum does not contain a half line, and it provides the first example of a Schr\"odinger operator in dimension $d \ge 2$ with any potential for which the spectrum is a generalized Cantor set.

One should nevertheless expect that for ``typical'' almost periodic Schr\"odinger operators in dimension $d \ge 2$, the spectrum contains a half line. On the other hand, the structure of the spectrum at low to medium energies is expected to be quite different. While there are absolutely no results in this direction, there is an extensive literature on almost periodic Schr\"odinger operators in one dimension and it is by now well understood that those operators tend to have Cantor spectrum; compare, for example, \cite{AS1981, DFL2017, E, M}. It is not unreasonable to expect that the spectra of ``typical'' almost periodic Schr\"odinger operators in dimension $d \ge 2$ do exhibit a Cantor structure at low and/or medium energies, in addition to the half line at high energies.

More generally one may ask whether there are any 
 Schr\"odinger operators in dimension $d \ge 2$ for which such a coexistence phenomenon, a Cantor structure and a half line in separate energy regimes, can be established. In fact, we were asked this precise question by Leonid Parnovski and the purpose of this paper is to provide an affirmative answer. Thus, our main result is the construction of the first examples of potentials in $\R^d$ exhibiting such a coexistence phenomenon.

Our examples will not be almost periodic in the classical (Bohr/Bochner) sense, but they belong to a class of uniformly recurrent potentials, which is defined as follows:

\begin{definition}\label{d.unifrec}
Let $d \ge 1$. A function $V : \R^d \to \R$ is called uniformly recurrent if for every $r > 0$, there exists $R > 0$ such that for every $x_0 \in \R^d$ and every $x_1 \in \R^d$, there exists $x_2 \in \R^d$ with $0 < |x_1-x_2| < R$ such that $V(x - x_0) = V(x - x_2)$ for every $x \in \R^d$ with $|x| < r$. The set of bounded uniformly recurrent $V : \R^d \to \R$ is denoted by $\mathcal{UR}_d$.
\end{definition}

\begin{theorem}\label{t.main}
For every $d \ge 2$, there are potentials $V \in \mathcal{UR}_d$ such that the spectrum of the associated Schr\"odinger operator in $L^2(\R^d)$ has a Cantor structure at low energies and possesses no gaps at high energies. More precisely, there are $E_0 < E_1$ such that $\sigma(H_V) \cap (-\infty, E_0]$ is perfect, nonempty, and nowhere dense and $\sigma(H_V) \cap [E_1, \infty) = [E_1, \infty)$.
\end{theorem}

\begin{remark}
(a) Uniform recurrence is related in spirit to classical almost periodicity. Both try to capture the fact that the sets of translates for which one observes suitable repetitions are relatively dense. The notions differ in the definition of a suitable repetition. Bohr almost periodicity requires global repetition up to a small $L^\infty$ error, while uniform recurrence requires exact repetition in a neighborhood of a given size. In particular, neither notion implies the other.
\\[2mm]
(b) Uniform recurrence is a standard notion in discrete geometry, and especially in the study of aperiodically ordered point sets, which in turn are of interest in the mathematical study of quasicrystals; see \cite{BG} and references therein.
\\[2mm]
(c) The potentials $V$ for which we prove these statements are completely explicit. They arise from a product construction based on the Fibonacci sequence, which may be expressed in terms of a discontinuous sampling function defined over an irrational rotation as follows:
\[
\omega_n =  \chi_{[1-\theta,1)}(n\theta \ \mathrm{mod} \ 1),
\quad \theta = \frac{\sqrt{5}-1}{2}.
\]
(d) It would of course be of interest to exhibit almost periodic $V : \R^d \to \R$, $d \ge 2$, for which the spectrum of the associated Schr\"odinger operator has both a Cantor component and a half line. We have been unable to establish the existence of such potentials. The obstacles which appear to make this task difficult are discussed in Section~\ref{s.?}; see Question~\ref{q.ap} and the discussion following it.
\end{remark}

\medskip

As already indicated above, our level of understanding is vastly different in the cases $d = 1$ and $d \ge 2$. Far more is known in the one-dimensional case, and in order to take advantage of that we will consider separable potentials $V : \R^d \to \R$, that is, potentials of the form
\begin{equation}\label{e.seppot}
V(x_1, \ldots, x_d) = V_1(x_1) + \cdots + V_d(x_d)
\end{equation}
with $V_j : \R \to \R$, $1 \le j \le d$. Thus, in addition to the operator $H_V$ in $L^2(\R^d)$ of main interest, we may also consider the operators $H_{V_j}$ in $L^2(\R)$, $1 \le j \le d$. It is well known that the spectra are related via the standard Minkowski sum of sets as follows:
\begin{equation}\label{e.seppotspec}
\sigma(H_V) = \sigma(H_{V_1}) + \cdots + \sigma(H_{V_d}) = \{E_1+\cdots + E_d : E_j \in \sigma(H_{V_j}), \ 1 \le j \le d\}.
\end{equation}
See, for example, \cite[Sections~II.4 and VIII.10]{RS1} and \cite[Theorem~II.10]{RS1}.

This suggests the strategy we will follow in the proof of Theorem~\ref{t.main}. Choose one-dimensional potentials $V_j$ such that  the spectra $\sigma(H_{V_j})$ are generalized Cantor sets. Then investigate the sum on the right-hand side of \eqref{e.seppotspec} and prove that there are $E_0, E_1$ such that it is Cantor below $E_0$ and a half line above $E_1$. It is well known that the study of sums of Cantor sets it difficult and can lead to a variety of outcomes. Such sums can indeed be Cantor sets or intervals, and certain mechanisms are known that will produce either outcome. The challenge will then be to choose the potentials $V_j$ in such a way that the spectra $\sigma(H_{V_j})$ have a structure that will allow us to apply both of these mechanisms simultaneously.

Once we have explained which features of Cantor sets we will need, we may explore these features in spectra of specific one-dimensional Schr\"odinger operators. Our one-dimensional potentials will be generated by the Fibonacci sequence. This is not a coincidence. This is currently (essentially) the only example for which the desired features can be established. In particular, no bona fide almost periodic potential in one dimension is currently known for which similar results can be proved. Thus, in Subsection~\ref{ss.fibonacci} we will discuss one-dimensional potentials generated by the Fibonacci substitution and recall some important concepts needed in the spectral analysis of the associated Schr\"odinger operators, namely the trace map, the Fricke-Vogt invariant, and the curve of initial conditions.

The remainder of the paper is organized as follows. In Section~\ref{s.absc}, we formulate and prove an abstract Bethe--Sommerfeld criterion for sums of extended Cantor sets. We single out this result, since it may be of independent interest outside our main application. In Section~\ref{s.back} we review some background about potentials generated by the Fibonacci substitution. In Section~\ref{s.hi}, we verify the conditions of the abstract Bethe--Sommerfeld criterion for the Fibonacci model and hence confirm the half-line portion of Theorem~\ref{t.main}. In Section~\ref{s.lo}, we study the low-energy region, showing that the spectrum can be made to be a Cantor set in suitable paramter regions. Finally, in Section~\ref{s.?} we conclude with some questions raised by the present work that we regard as challenging and interesting.

\begin{remark}
(a) For notational simplicity, we give the proof of Theorem~\ref{t.main} in the case $d = 2$. The proof in the case $d > 2$ is completely analogous. In fact, no changes are necessary in the choice of the underlying 1D potential to get the half-line portion of the spectrum, and only a minor change is necessary to get the Cantor portion of the spectrum (namely, one needs to replace the number $1/2$ in \eqref{e.smengoal2} by $1/d$).
\\[2mm]
(b) It is interesting to compare the structure of the spectrum of Schr\"odinger operators given by Theorem \ref{t.main} with the structure of Lagrange and Markov spectra, the sets that appear naturally in the theory of Diophantine approximation; see \cite{Math} and references therein.
\end{remark}

\bigskip

\noindent\textbf{Acknowledgment.} We are grateful to Leonid Parnovski for raising the question of whether examples with spectra exhibiting the topological structure discussed in this paper can be found and for his interest in this work.

\section{Abstract Bethe--Sommerfeld Criterion}\label{s.absc}

Let us begin by recalling some terminology and results that will be useful. Let $C \subset \mathbb{R}$ be a Cantor set (i.e., $C$ is bounded, perfect, and nowhere dense) and denote by $I$ its convex hull. Any connected component of $I\backslash C$ is called a \emph{gap} of $C$. A \emph{presentation} of $C$ is given by an ordering $\mathcal{U} = \{U_n\}_{n \ge 1}$ of the gaps of $C$. If $u \in C$ is a boundary point of a gap $U$ of $C$, we denote by $K$ the connected component of $I\backslash (U_1\cup U_2\cup\ldots \cup U_n)$ (with $n$ chosen so that $U_n = U$) that contains $u$ and write
$$
\tau(C, \mathcal{U}, u)=\frac{|K|}{|U|}.
$$
The \emph{thickness} $\tau(C)$ of $C$ is given by
$$
\tau(C) = \sup_{\mathcal{U}} \inf_{u} \tau(C, \mathcal{U}, u), 
$$

The following consequence of the Newhouse Gap Lemma \cite{N70, N79} is stated as \cite[Lemma~6.2]{DG2011} and proved there:

\begin{lemma}\label{l.sumofcs}
Suppose $C, K\subset \mathbb{R}$ are Cantor sets with $\tau(C)\cdot \tau(K)>1$. Assume also that the size of the largest gap of $C$ is not greater than the diameter of $K$, and the size of the largest gap
of $K$ is not greater than the diameter of $C$. Then,
$$
C + K = [\min C + \min K, \max C + \max K].
$$
\end{lemma}

\begin{remark}\label{r.cplusc}
A particular consequence of Lemma~\ref{l.sumofcs} is the following: if $C \subset \mathbb{R}$ is a Cantor set with $\tau(C) > 1$, then
$$
C + C = [2 \min C, 2 \max C].
$$
\end{remark}

\begin{lemma}[Abstract Bethe-Sommerfeld Conditions]\label{l.ABSC}
Suppose $K\subset \R$ is a closed set that has the following properties:
\begin{enumerate}
\item $K$ contains Cantor sets $K_1, K_2, \ldots$ with disjoint convex hulls $I_1, I_2, \ldots$, so that $I_{n+1}$ lies to the right of $I_{n}$ for every $n \in \mathbb{N}$.
\item For some $\varepsilon > 0$, we have $\tau(K_n) > 1 + \varepsilon$ for all $n \in \mathbb{N}$.
\item For some $A > a > 0$, we have $2A \ge |I_n| \ge A$ and $\mathrm{dist} (I_n, I_{n+1}) \le a$ for all $n \in \mathbb{N}$.
\end{enumerate}
Let $\widetilde K=f(K)$, where $f:\mathbb{R} \to \mathbb{R}$ is given by $f(x)=x^2$.  Then $\widetilde K+\widetilde K$ contains a half-line.
\end{lemma}

\begin{proof}
Denote $\widetilde K_n=f(K_n)$. Throughout the argument, we assume $n$ is large enough that $K_n \subseteq [0,\infty)$. Let us show that for all sufficiently large $n\in \mathbb{N}$, we have $\tau(\widetilde K_n)\ge 1+\frac{\varepsilon}{2}$. Indeed, for any two intervals $B=[u,v]$, $U=[v, w]$ from $I_n$ such that
$$
\frac{|B|}{|U|}=\frac{v-u}{w-v}\ge 1+\varepsilon,
$$
we have
\begin{align*}
\frac{|f(B)|}{|f(U)|} & = \frac{v^2 - u^2}{w^2 - v^2} \\
& = \frac{v - u}{w - v} \cdot \frac{v + u}{w + v} \\
& > (1 + \varepsilon)\left( \frac{2v - 2A}{2v+2A}\right) \\
& = (1 + \varepsilon) \left( 1 - \frac{2A}{v+A} \right) \\
& > 1 + \frac{\varepsilon}{2},
\end{align*}
for $v$ sufficiently large (i.e.\ $n$ sufficiently large).

 A similar estimate holds in the case of intervals $B=[u,v]$, $U=[w, u]$ (i.e., in the language of the Newhouse Gap Lemma, when the bridge is to the right of the gap). Therefore, for all sufficiently large $n$, $\tau(\widetilde K_n) \ge 1+ \frac{\varepsilon}{2}$, and thus the sum $\widetilde K_n+\widetilde K_{n}$ is an interval by Remark~\ref{r.cplusc}.

Next, let us show that for all sufficiently large $n$, each of the Cantor sets $\widetilde K_n$ and  $\widetilde K_{n+1}$ has diameter larger than the largest gap of the other. In that case Lemma~\ref{l.sumofcs} will imply that $\widetilde K_n+\widetilde K_{n+1}$ is an interval. Indeed, suppose that $I_n=[x,y]$ and $I_{n+1}=[z,t]$. By our assumptions, we have
\begin{equation} \label{eq:ABSC:basicIneqs}
A \le y-x \le 2A, \quad z-y \le a, \quad A \le t-z \le 2A.
\end{equation}
Since we already know that $\tau(\widetilde K_{n+1})>1$, the largest gap of $\widetilde K_{n+1}$ is not greater than $\frac{1}{3}(t^2-z^2)$, and the diameter of $\widetilde K_n$ is equal to $y^2-x^2$. Using \eqref{eq:ABSC:basicIneqs} three times, we have
\begin{align*}
y^2 - x^2 & = (y - x)(y + x)\\
& \ge A(x + y) \\
& \ge A(2y-2A) \\
& \ge A(2(z - a) - 2A) \\
& = 2Az - (2Aa + 2A^2).
\end{align*}
Similarly, using \eqref{eq:ABSC:basicIneqs} twice more yields
\begin{align*}
\frac{1}{3}(t^2 - z^2) & = \frac{1}{3}(t - z)(t + z) \\
& \le \frac{1}{3} 2A(t + z) \\
& \le \frac{2}{3} A(2z + 2A) \\
& = \frac{4}{3} Az + \frac{4}{3} A^2.
\end{align*}
Therefore, for all sufficiently large values of $n$, the diameter of $\widetilde K_n$ is greater than the largest gap of $\widetilde K_{n+1}$. Similarly one can show that for all sufficiently large values of $n$, the diameter of $\widetilde K_{n+1}$ is greater than the largest gap of $\widetilde K_{n}$.

We therefore know that the sets $J_n := \widetilde K_n + \widetilde K_{n}$ and $J_n' := \widetilde K_n + \widetilde K_{n+1}$ are intervals for large $n$. Let us show that they cover a half line. It is enough to check that $J_n$ is not disjoint from $J_n'$, and $J_n'$ is not disjoint from $J_{n+1}$.

As before, denote $I_n = [x, y]$, $I_{n+1} = [z,t]$. It follows from our discussion above that $J_n = [2x^2, 2y^2]$, $J_{n+1} = [2z^2,2t^2]$, and $J_n'=[x^2 + z^2, y^2 + t^2]$. To show that $J_n$ is not disjoint from $J_n'$ we need to check that $2y^2\ge x^2+z^2$. Using \eqref{eq:ABSC:basicIneqs} twice in each term, we have
\begin{align*}
2 y^2 - (x^2 + z^2)
& = (y-x)(y+x) + (y-z)(y+z) \\
& \ge  A(y + x) - a(y + z) \\
& \ge A(2y - 2A) - a(2y + a) \\
& = 2y(A - a) - 2A^2 - a^2 \\
& \ge 0
\end{align*}
if $n$ is sufficiently large.

To show that $J_n'$ is not disjoint from $J_{n+1}$ we need to check that $y^2 + t^2 \ge 2z^2$. We have
\begin{align*}
y^2 + t^2-2z^2 & = (t^2-z^2)+(y^2-z^2) \\
& = (t-z)(t+z)+(y-z)(y+z).
\end{align*}
Using $t-z = |I_{n+1}| \geq A$ and $y \leq z \leq t$ to estimate the first term and  $z-y \leq a$ (twice) to estimate the second term, we arrive at
\begin{align*}y^2 + t^2-2z^2
& \ge 2Ay - a(2y+a) \\
& \ge 2(A-a)y - a^2 \\
& \ge 0
\end{align*}
if $n$ is sufficiently large.

It follows that for $n_0$ large enough, the set
$$
\bigcup_{n \ge n_0} \left( J_n \cup J_n' \right),
$$
which is contained in $\widetilde K+\widetilde K$, is a half line.
\end{proof}

\section{Background} \label{s.back}

\subsection{Trace Map and Fricke--Vogt Invariant}\label{ss.fibonacci}

Let us recall the dynamical setup that is canonically associated with any one-dimensional Schr\"odinger operator whose potential is generated by the Fibonacci substitution. We refer the reader to \cite{DEG2015, DG2011} for more details and background.

The \textit{trace map} is given by
$$
T : \Bbb{R}^3 \to \Bbb{R}^3, \; T(x,y,z)=(2xy-z,x,y).
$$
Note that $T$ is invertible with $T^{-1}(x,y,z) = (y, z, 2yz - x)$.

The \textit{Fricke-Vogt invariant}
\begin{equation*}
G(x,y,z) = x^2+y^2+z^2-2xyz-1
\end{equation*}
obeys
\begin{equation}\label{e.fvinvariant2}
G \circ T = G,
\end{equation}
and hence $T$ preserves the family of cubic surfaces
$$
S_I = \left\{ (x,y,z) \in \R^3 : G(x,y,z) = I \right\}.
$$
The surface $S_0$ is called the \textit{Cayley cubic}.

It is natural to consider the restriction $T_{I}$ of the trace map $T$ to the invariant surface $S_I$. That is, $T_{I} : S_I \to S_I$, $T_{I} = T|_{S_I}$. Denote by $\Omega_{I}$ the set of points in $S_I$ whose full two-sided orbit under $T_{I}$ is bounded.

Let us recall that an invariant closed set $\Lambda$ of a diffeomorphism $f : M \to M$ is \textit{hyperbolic} if there exists a splitting of the tangent space $T_x M = E^u_x \oplus E^u_x$ at every point $x \in \Lambda$ such that this splitting is invariant
under $Df$, the differential $Df$ exponentially contracts vectors from the stable subspaces $\{E^s_x\}$, and the differential of the inverse, $Df^{-1}$, exponentially contracts vectors from the unstable subspaces $\{E^u_x\}$. A hyperbolic set $\Lambda$ of a
diffeomorphism $f : M \to M$ is \textit{locally maximal} if there exists a neighborhood $U$ of $\Lambda$ such that
$$
\Lambda=\bigcap_{n\in\Bbb{Z}}f^n(U).
$$

It was shown in \cite{Can, Cas, DG2009} that for every $I > 0$, the set $\Omega_{I}$ is a locally maximal hyperbolic set of $T_{I}:S_I \to S_I$, which is homeomorphic to a Cantor set.

For the discussion of the high-energy regime later in this paper, special attention needs to be paid to small non-negative values of $I$, which may then be viewed as small perturbations of the case $I = 0$. For this reason, let us briefly discuss the latter special case.

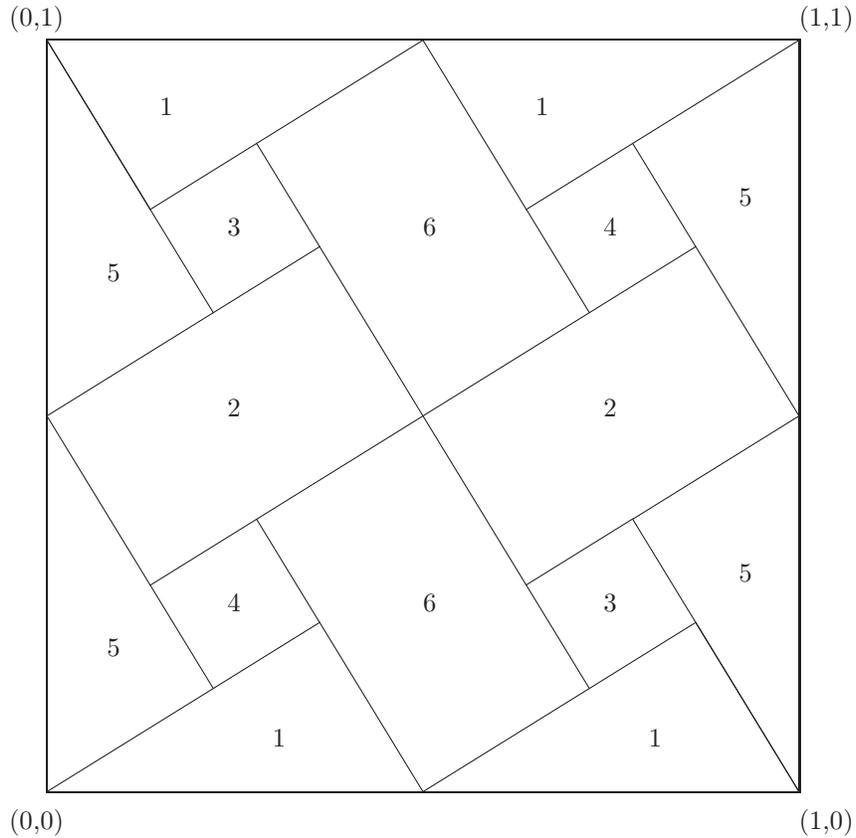
\begin{figure}[htb]
\setlength{\unitlength}{1mm}
\begin{picture}(120,120)

\put(10,10){\framebox(100,100)}

\put(5,5){(0,0)}

\put(110,5){(1,0)}

\put(5,112){(0,1)}

\put(110,112){(1,1)}

\put(10,10){\line(161,100){36.2}}

\put(60,10){\line(-61,100){13.8}}

\put(40,16){$1$}

\put(60,10){\line(161,100){36.2}}

\put(110,10){\line(-61,100){13.8}}

\put(90,16){$1$}

\put(60,110){\line(-161,-100){36.2}}

\put(10,110){\line(61,-100){13.8}}

\put(25,100){$1$}

\put(110,110){\line(-161,-100){36.2}}

\put(60,110){\line(61,-100){13.8}}

\put(75,100){$1$}

\put(10,60){\line(61,-100){22.1}}

\put(18,28){$5$}

\put(10,60){\line(161,100){22.1}}

\put(10,110){\line(61,-100){22.1}}

\put(18,78){$5$}

\put(110,60){\line(-61,100){22.1}}

\put(102,38){$5$}

\put(110,60){\line(-161,-100){22.1}}

\put(110,10){\line(-61,100){22.1}}

\put(102,88){$5$}

\put(60,60){\line(-161,-100){36.3}}

\put(60,60){\line(161,100){36.3}}

\put(60,60){\line(-61,100){22.1}}

\put(60,60){\line(61,-100){22.1}}

\put(88,46.35){\line(-161,-100){14.3}}

\put(32,73.65){\line(161,100){14.3}}

\put(37.88,46.28){\line(61,-100){8.4}}

\put(82.12,73.72){\line(-61,100){8.4}}

\put(34,34){$4$}

\put(84,84){$4$}

\put(34,84){$3$}

\put(84,34){$3$}

\put(34,60){$2$}

\put(84,60){$2$}

\put(60,34){$6$}

\put(60,84){$6$}

\end{picture}
\caption{The Markov partition for the map $\mathcal{A}$.}\label{fig.Casdagli-Markov}
\end{figure}

Denote by $\mathbb{S}$ the part of the Cayley cubic $S_0$ inside the cube $\{ |x|\le 1, |y|\le 1, |z|\le 1\}$. The surface $\mathbb{S}$ is homeomorphic to $S^2$, invariant, smooth everywhere except at the four points $P_1=(1,1,1)$, $P_2=(-1,-1,1)$, $P_3=(1,-1,-1)$, and $P_4=(-1,1,-1)$, where $\mathbb{S}$ has conic singularities, and the trace map $T$ restricted to $\mathbb{S}$ is a factor of the hyperbolic automorphism of $\T^2 = \R^2 / \Z^2$ given by
$$
\mathcal{A}(\theta, \varphi)=(\theta + \varphi, \theta)\ (\text{\rm mod}\ 1).
$$
The semiconjugacy is given by the map
$$
F: (\theta, \varphi) \mapsto (\cos 2\pi(\theta + \varphi), \cos 2\pi \theta, \cos 2\pi \varphi).
$$
The map $\mathcal{A}$ is hyperbolic, and is given by the matrix $A=\begin{pmatrix} 1 & 1 \\ 1 & 0 \end{pmatrix}$, which has
eigenvalues
$$
\mu=\frac{1+\sqrt{5}}{2}\ \ \text{\rm and} \ \ \ -\mu^{-1}=\frac{1-\sqrt{5}}{2}.
$$
The Markov partition for the map $\mathcal{A}:\mathbb{T}^2\to \mathbb{T}^2$ shown in Figure~\ref{fig.Casdagli-Markov} has already appeared in \cite{Cas, DG2011}; for more details on Markov partitions for two-dimensional hyperbolic maps see \cite[Appendix 2]{PT}). Its image under the map $F:\mathbb{T}^2\to \mathbb{S}$ is a Markov partition for the pseudo-Anosov map $T_0 : \mathbb{S}\to \mathbb{S}$.

So far we have discussed the dynamical setup which is tied to the trace map. Since all ``second-order'' one-dimensional operators (more precisely, operators for which there exists an $\mathrm{SL}(2,\R)$ transfer matrix formalism) generated by the Fibonacci substitution lead to the same trace map, this dynamical setup is relevant to all of them. The specific choice of the operator in question, however, will affect the associated curve of initial conditions, which needs to be studied relative to the fixed dynamical setup. Let us discuss the operator and its associated curve of initial conditions relevant to this paper.

\subsection{Schr\"odinger Operators Associated with the Fibonacci Subshift} \label{s:subshift}
Let us introduce the second-order operators that we will study. See also \cite{DFG2014,FM}. First, we recall the notion of concatenation of real-valued functions defined on intervals. Given $\ell_n \in \R_+ \eqdef (0,\infty)$ and $f_n:[0,\ell_n) \to \R$ for each $n \in \Z$, define the \emph{concatenation} of the sequence $\{f_n\}_{n \in \Z}$ as follows. Put
\begin{equation} \label{eq:concatEndptDef}
t_n
\eqdef
\begin{cases}
\sum_{j=0}^{n-1} \ell_j & n \geq 1 \\
0 & n = 0 \\
-\sum_{j=n}^{-1} \ell_j & n \leq -1,
\end{cases}
\end{equation}
denote $J_n = [t_{n}, t_{n+1})$, $J = \bigcup_n J$ and define $f:J \to \R$ by
\begin{equation}\label{eq:concatDef}
f(x)
=
f_n(x - t_{n}),\quad
\text{ for each } x \in J_n.
\end{equation}
In the present work, we will have $\sum_{n<0} \ell_n = \sum_{n > 0} \ell_n = \infty$ so that $J = \R$.
Using a box to denote the position of the origin, we denote the concatenation by
\begin{equation} \label{eq:finiteConcatDef}
f
=
\left( \cdots \; | \; f_{-2} \; | \; f_{-1} \; | \; \boxed{f_0} \; | \; f_1 \; | \; f_2 \; | \;  \cdots  \right).
\end{equation}

Let $\alphabet$ be a finite set, called the \emph{alphabet}. Equip $\alphabet$ with the discrete topology and endow $\alphabet^\Z$ with the corresponding product topology. The \emph{left shift}
\[
[T \omega](n) \eqdef \omega(n+1),
\quad \omega \in \alphabet^{\Z}, \; n \in \Z,
\]
defines a homeomorphism from $\alphabet^{\Z}$ to itself. A subset $\Omega \subseteq \alphabet^\Z$ is called \emph{$T$-invariant} if $T^{-1}(\Omega) = \Omega$. Any compact $T$-invariant subset of $\alphabet^\Z$ is called a \emph{subshift}.

We can use the concatenation construction above to associate potentials (and hence Schr\"odinger operators) to elements of subshifts as follows. For each $\alpha \in \alphabet$, we pick $\ell_\alpha > 0$ and a real-valued function $f_\alpha \in L^2[0,\ell_\alpha)$. Then, for any $\omega \in \alphabet^\Z$, we define the action of the continuum Schr\"odinger operator $H_\omega$ in $L^2(\R)$ by
\begin{equation}\label{eq:Schrod}
H_\omega
=
- \frac{\dd^2}{\dd x^2}+ V_\omega \cdot,
\end{equation}
where the potential $V_\omega$ is given by
\begin{equation} \label{eq:VomegaDef}
V_{\omega}
=
V_{\omega,\{f_\alpha\}}
\eqdef
\left( \cdots \; f_{\omega_{-2}} \; | \; f_{\omega_{-1}} \; | \; \boxed{f_{\omega_0}} \; | \; f_{\omega_1} \; | \; f_{\omega_2} \; \cdots \right).
\end{equation}
These potentials belong to $L^2_\mathrm{\loc, \mathrm{unif}}(\R)$ and hence each $H_\omega$ defines a self-adjoint operator on a dense subspace $L^2(\R)$ in a canonical fashion.

 In the present manuscript, the alphabet contains two symbols, $\alphabet \eqdef \{ \a,\b \}$. The Fibonacci substitution is the map
\[
S(\a) = \a\b, \; S(\b) = \a.
\]
This map extends by concatenation to $\alphabet^*$, the free monoid over $\alphabet$ (i.e.\ the set of finite words over $\alphabet$), as well as to $\alphabet^{\N}$, the collection of (one-sided) infinite words over $\alphabet$. There exists a unique element
\[
u
=
\a\b\a\a\b\a\b\a\ldots \in \alphabet^{\N}
\]
with the property that $u = S(u)$. It is straightforward to verify that for $n \in \N$, $S^n(\a)$ is a prefix of $S^{n+1}(\a)$. Thus, one obtains $u$ as the limit (in the product topology on $\alphabet^\mathbb{N}$) of the sequence of finite words $\{ S^n(\a) \}_{n \in \N}$. With this setup, the Fibonacci subshift is defined to be the collection of two-sided infinite words with the same local factor structure as $u$, that is,
\[
\Omega
\eqdef
\{ \omega \in \alphabet^\Z : \text{every finite subword of $\omega$ is also a subword of } u \}.
\]
Given $\ell_\a, \ell_\b \in \R_+$ and real-valued functions $f_\a \in L^2[0,\ell_\a)$, $f_\b \in L^2[0,\ell_\b)$, we consider the family of continuum Schr\"odinger operators $\{H_\omega\}_{\omega \in \Omega}$ defined by \eqref{eq:Schrod} and \eqref{eq:VomegaDef}. Since $(\Omega,T)$ is a minimal dynamical system, one can verify that there is a uniform closed set $\Sigma = \Sigma(f_\a,f_\b) \subset \R$ with the property that
\begin{equation}\label{e.sigmaisconstant}
\sigma(H_\omega)
=
\Sigma
\text{ for every }
\omega \in \Omega.
\end{equation}

Of course, one can choose $f_\a$ and $f_\b$ in such a way that every $V_\omega$ is a periodic potential (notice that as soon as $V_{\omega_0}$ is periodic for a single $\omega_0 \in \Omega$, then every $V_\omega$ is periodic by minimality). The main result of \cite{DFG2014} is that this is the only possible obstruction to Cantor spectrum. Thus, we adopt the following assumption throughout the paper:

\begin{assumption} \label{assumption}
The potential pieces $f_\a$ and $f_\b$ are chosen so that $V_\omega$ is aperiodic for one $\omega \in \Omega$ (hence for every $\omega \in \Omega$ by minimality).
\end{assumption}

\begin{theorem}[Damanik--Fillman--Gorodetski \cite{DFG2014}]
If $f_\a$ and $f_\b$ satisfy Assumption~\ref{assumption}, then $\Sigma(f_\a,f_\b)$ is a Cantor set of zero Lebesgue measure.
\end{theorem}

\begin{remark}
In \cite{DFG2014}, the authors also assumed a condition on $\{f_\alpha:\alpha\in \alphabet\}$ that they called \emph{irreducibility}. This condition is defined so that the potentials satisfy the \emph{simple finite decomposition property} (SFDP) from \cite{KLS2011}. However, since our alphabet only has two letters, SFDP follows from aperiodicity and \cite[Proposition~3.5]{KLS2011}.
\end{remark}

In the sequel, it will be convenient to introduce an additional symbol: ``$\a\b$'', and then to define $\ell_\ab = \ell_\a+\ell_\b$ and
\[
f_\ab(x)
=
\begin{cases}
f_\a(x) & 0 \leq x < \ell_\a \\
f_\b(x-\ell_\a) & \ell_\a \leq x < \ell_\ab.
\end{cases}
\]

Given $\alpha \in \widetilde{\alphabet} \eqdef \alphabet \cup\{\ab\}$, consider the solutions of the differential equation $-u''(x) + f_\alpha(x) u(x) = E u(x)$ for  $E \in \R$. Denote the solution obeying $u(0) = 0$, $u'(0) = 1$ (resp., $u(0) = 1$, $u'(0) = 0$) by $u_{\alpha,D}(\cdot,E)$ (resp., $u_{\alpha,N}(\cdot,E)$). Then, we set
\begin{align*}
M(\alpha,E) & = \begin{pmatrix} u_{\alpha,N}(\ell_\alpha,E) & u_{\alpha,D}(\ell_\alpha,E) \\
 u_{\alpha,N}'(\ell_\alpha,E) & u_{\alpha,D}'(\ell_\alpha,E) \end{pmatrix},
\end{align*}
and
\begin{align}
x_{-1}(E) & = \frac12 \tr \left( M(\b,E) \right), \label{e.x-1} \\
x_0(E) & = \frac12 \tr \left( M(\a,E) \right), \label{e.x0} \\
x_1(E) & = \frac12 \tr \left( M(\ab,E)  \right) . \label{e.x1}
\end{align}
It is straightforward to check that $M(\ab,\cdot) = M(\b,\cdot)M(\a,\cdot)$. The map $\gamma : E \mapsto (x_1(E), x_0(E), x_{-1}(E))$ will be called the \textit{curve of initial conditions}. Note that by \eqref{e.fvinvariant2}, all the points $T^n(x_1(E), x_0(E), x_{-1}(E))$ lie on the surface $S_{I(E)}$, where (with some abuse of notation) we set
$$
I(E) = I(\gamma(E)).
$$
It was shown in \cite{DFG2014} that $I(E) \ge 0$ for every $E \in \Sigma$.

For special choices of the local potential pieces, it is possible to compute the Fricke-Vogt invariant $I(E)$ explicitly. For example, as in \cite{DFG2014} let us consider the case $\ell_\a=\ell_\b = 1$, $f_\a = \lambda \cdot \chi_{[0,1)}$ and $f_\b = 0 \cdot \chi_{[0,1)}$, where $\lambda > 0$ and let us recall the formulas obtained there. Clearly, the resulting potentials are aperiodic. One readily computes the traces as follows:
\begin{align*}
x_{-1}(E) & = \cos \sqrt{E}, \\
x_0(E) & = \cos \sqrt{E - \lambda}, \\
x_1(E) & = \cos \sqrt{E} \cos \sqrt{E-\lambda} -\frac{1}{2}\left(\sqrt{\frac{E}{E-\lambda}}+\sqrt{\frac{E-\lambda}{E}}\right)\sin \sqrt{E} \sin \sqrt{E-\lambda},
\end{align*}
and therefore
\begin{equation}\label{e.1}
I(E) = \frac{1}{4}\frac{\lambda^2}{E(E-\lambda)}\sin^2 \sqrt{E} \sin^2 \sqrt{E-\lambda}.
\end{equation}
These expressions all define entire functions of $E \in \C$. For example, for $E < 0$, one has $x_{-1}(E) = \cosh\sqrt{|E|}$.



Let us now state our main 1D results which we will leverage to prove Theorem~\ref{t.main}. Let $\ell_\a=\ell_\b=1$ and $f_\a = 0\cdot \chi_{[0,1)}$ and $f_\b = \lambda\cdot \chi_{[0,1)}$ be as above, and write $\Sigma_\lambda = \Sigma( f_\a, f_\b)$.

\begin{theorem} \label{t:hi}
For any $\lambda \in \R$, $\Sigma_\lambda + \Sigma_\lambda$ contains a half-line. More precisely, there exists $E_1 = E_1(\lambda)$ such that $\Sigma_\lambda + \Sigma_\lambda \supseteq [E_1,\infty)$.
\end{theorem}

\begin{theorem} \label{t:lo}
There exists $\lambda_0$ such that, for any $\lambda \geq \lambda_0$, there is $E_0$ such that $(-\infty,E_0] \cap (\Sigma_\lambda+\Sigma_\lambda)$ is a nonempty Cantor set of zero Lebesgue measure.
\end{theorem}

\begin{proof}[Proof of Theorem~\ref{t.main}]
Note first that it is well known that each $\omega \in \Omega$ is uniformly recurrent (where the definition of uniform recurrence for maps $\omega : \Z \to \R$ is analogous to the one given in Definition~\ref{d.unifrec} for maps $V : \R \to \R$; namely, for every length $m \in \N$, there is a window size $M \in \N$ such that each word of length $M$ appearing in $\omega$ contains every word of length $m$ that appears in $\omega$).

From this and the way the potentials $V_\omega$ are constructed, it is then easily seen that each $V_\omega$ is uniformly recurrent in the sense of Definition~\ref{d.unifrec}. The resulting $2$-dimensional potential arising via \eqref{e.seppot} with $V_1 = V_2 = V_\omega$ then inherits the uniform recurrence property from $V_\omega$. (Recall that we focus on the case $d = 2$ in this proof and leave the details of the extension to the case $d \ge 3$ to the reader.)

The statements about the spectrum are immediate from Theorem~\ref{t:hi}, Theorem~\ref{t:lo}, and Equation~\eqref{e.seppotspec}.
\end{proof}

\section{The High-Energy Regime} \label{s.hi}

Our goal in this section is to prove Theorem~\ref{t:hi}, which states that in the high-energy regime, the spectrum of the separable 2D Schr\"odinger operators we consider contains no gaps. In this proof the existing (extensive) work on the Fibonacci Hamiltonian will be helpful, but not sufficient. We will actually have to delve into the technical details of the theory that has been developed and establish new results. For this reason, this section is somewhat demanding and while we attempt to provide a presentation that is as self-contained as possible, it may be necessary to consult some of the earlier papers (e.g., \cite{DFG2014, DG2009, DG2011, DGY2016, FTY, Y}) for further background and motivation.

The curve of initial conditions for the degenerate case in which the potential vanishes identically (i.e., for the free Laplacian) contains the curve $\Pi = \{ (\cos (2t), \cos t, \cos t) \}$, $t \ge 0$. This curve corresponds to the line segment in Figure~\ref{fig.Casdagli-Markov} that connects the point $(0,0)$ with the point $(0.5, 0.5)$. In particular, this curve is transversal to both the stable and unstable foliations of the trace map $T_0$ on $\S_0$. The curve of initial conditions $\Gamma(E)$ accumulates to the curve $\Pi$ as $E$ becomes large, if one sets $t = \sqrt{E}$. We would like to consider only the parts of $\Gamma(E)$ that are bounded away from suitably small neighborhoods of the singularities, and to then show that on each connected piece of the curve, the set of points that have bounded positive semi-orbits (i.e., those that correspond to energies in the spectrum) form a Cantor set of large thickness. In order to formalize this construction we need to discuss the dynamics of the trace map $T$ and to introduce some notation.

First consider the dynamics of $T$ in a neighborhood of the singularity $P_1 = (1,1,1)$. Due to the symmetries of the trace map this will also provide information on the dynamics near the other singularities. Take $r_0 > 0$ small and let $O_{r_0}(P_1)$ be an $r_0$-neighborhood of the point $P_1 = (1, 1, 1)$ in $\mathbb{R}^3$. Let us consider the set $\Per_2(T)$ of periodic points of $T$ of period 2.

%
%
%

By a direct calculation one can show the following:

\begin{lemma}\label{l.periodtwo} We have
$$
\Per_2(T) = \left\{ (x,y,z) : x \in (-\infty, \tfrac{1}{2}) \cup (\tfrac{1}{2}, \infty), \ y = \frac{x}{2x-1}, \ z = x \right\}.
$$
\end{lemma}

Notice that in a neighborhood of $P_1$, the intersection $\mathfrak{I} \equiv \Per_2(T) \cap O_{r_0}(P_1)$ is a smooth curve that is a normally hyperbolic with respect to $T$ (see, e.g., Appendix~1 in \cite{PT} for the formal definition of normal hyperbolicity). Therefore, the local center-stable manifold  $W_{\loc}^{\cs}(I)$ and the local center-unstable manifold $W^{\cu}_{\loc}(I)$ defined by
$$
W_{\loc}^{\cs}(I) = \left\{ p \in O_{r_0}(P_1) : T^n(p) \in O_{r_0}(P_1) \ \text{\rm for all}\ n \in \mathbb{N} \right\},
$$
$$
W_{\loc}^{\cu}(I) = \left\{ p \in O_{r_0}(P_1) : T^{-n}(p) \in O_{r_0}(P_1) \ \text{\rm for all}\ n \in \mathbb{N} \right\}
$$
are smooth two-dimensional surfaces. Also, the local strong stable manifold $W^{\sss}_{\loc}(P_1)$ and the local strong unstable manifold $W^{\uu}_{\loc}(P_1)$ of the fixed point $P_1$, defined by
$$
W^{\sss}_{\loc}(P_1) = \left\{ p \in W^{\cs}_{\loc}(I) : T^n(p) \to P_1 \ \text{\rm as}\ n \to \infty \right\},
$$
$$
W^{\uu}_{\loc}(P_1) = \left\{ p \in W^{\cu}_{\loc}(I) : T^{-n}(p) \to P_1 \ \text{\rm as}\ n \to \infty \right\},
$$
are smooth curves.

The Markov partition for the pseudo-Anosov map $T : \mathbb{S} \to \mathbb{S}$ can be extended to a Markov partition for the map $T : S_I \to S_I$. Namely, there are four singular points $P_1 = (1,1,1)$, $P_2 = (-1,-1,1)$, $P_3 = (1,-1,-1)$, and $P_4 = (-1,1,-1)$ of $\mathbb{S}$. The point $P_1$ is a fixed point of $T$, and the points $P_2$, $P_3$, $P_4$ form a periodic orbit of period $3$. For small $I > 0$, on the surface $S_I$ near $P_1$ there is a hyperbolic orbit of the map $T_I = T|_{S_I}$ of period $2$, and near the orbit $\{P_2, P_3, P_4\}$ there is a hyperbolic periodic orbit of period $6$. Pieces of stable and unstable manifolds of these $8$ periodic points form a Markov partition for $T_I : S_I \to S_I$. For $I > 0$, the elements of this Markov partition are disjoint. Let us denote these six rectangles (the elements of the Markov partition) by $R_I^1, R_I^2, \ldots, R_I^6$. Let us also denote $R_I = \bigcup_{i=1}^6 R_I^i$ and $R = \bigcup_{0 \le I} R_I$. Finally, let us denote $R^i = \bigcup_{0 \le I} R_I^i$, $i = 1, 2, \ldots, 6$.

Denote by $\Omega^+$ the set of points in $\{I\ge 0\}\subset \mathbb{R}^3$ with bounded positive semiorbits.

The set of fixed points of $T^2$ in a neighborhood of $P_1$ of size $r_0 > 0$ is a smooth curve $\text{\rm Fix}(T^2, O_{r_0}(P_1)) = \Per_2(T) \cap O_{r_0}(P_1)$; see Lemma~\ref{l.periodtwo} above.  
Each of the fixed points has one of the eigenvalues equal to $1$, one greater than $1$, and one smaller than $1$ in absolute value.
Therefore the curve $\text{\rm Fix}(T^2, O_{r_0}(P_1))$ is a normally hyperbolic manifold, and its stable set $W^{\mathrm s}(\text{\rm Fix}(T^2, O_{r_0}(P_1)))$ is a smooth two dimensional surface; see \cite{HPS}. The strong stable manifolds form a $C^1$-foliation of $W^{\mathrm s}(\text{\rm Fix}(T^2, O_{r_0}(P_1)))$; see \cite[Theorem B]{PSW}.  It is convenient to consider $T^6$ since in this case each of the eight periodic points that were born from the singularities becomes a fixed point. Due to the symmetries of the trace map, the dynamics of $T^6$ is the same in a neighborhood of each of the singularities $P_1, P_2, P_3,$ and $P_4$.

In order to specify how to choose the pieces of the curve $\Gamma(t)$, let us consider an image $T^{6m}(\Pi)$. This is a curve that connects $P_1$ with $P_3$, and each connected component of the intersection of this curve with one of the elements of the Markov partition is a curve that connects one stable boundary with another. Let us take the curve obtained from  $T^{6m}(\Pi)$ by cutting the first and the last pieces with respect to this splitting. If $m$ is large, the preimage of this curve is $\Pi$ with small pieces near the singularities removed; let us denote it by $\widetilde \Pi$. The curve $\widetilde \Pi$ is parameterized by $t \in [2\pi n + \alpha, (2n+1) \pi - \beta]$ (or by $t\in [(2n+1) \pi + \beta, 2\pi n - \alpha]$) for some small $\alpha, \beta > 0$. In the curve $\Gamma(t)$ one can find a sequence of pieces such that the gaps between them are small, and the $T^{6m}$-images of these pieces connect the points of the boundaries of the same $R^i$ as $T^{6m}(\widetilde \Pi)$. Denote those pieces by $\Gamma_n$.

\begin{prop}\label{p.thick}
Let $\{\Gamma_n\}$ be a sequence of curves in $\mathbb{R}^3$ parameterized by $t \in [2\pi n + \alpha_n, (2n+1) \pi - \beta_n] \equiv J_n$, $\alpha_n \to \alpha, \beta_n \to \beta$, such that
\begin{enumerate}

\item $\Gamma_n$ converges to $\widetilde \Pi$ in the $C^2$ topology;

\item $I(\Gamma_n(t)) > 0$ for all $t \in J_n$;

\item for some uniform $C>0$, we have $\left| \frac{dI(\Gamma_n(t))}{dt}\right| \le C I(\Gamma_n(t))$;

\item the images $T^k(\Gamma_n)$ intersect exactly the same sets $R^i$ (and in the same order) as $T^k(\widetilde \Pi)$ for all $k\ge 0$.

\end{enumerate}
Then the set $K_n = \{ t \in J_n : \Gamma_n(t) \in \Omega^+ \}$ is a Cantor set for all large values of $n$, and $\tau(K_n) \to \infty$ as $n \to \infty$.
\end{prop}

The proof of Proposition~\ref{p.thick} goes along similar lines as the proof of \cite[Theorem~1.2]{DG2011}, which asserts that the spectrum of the discrete Fibonacci Hamiltonian is a Cantor set of large thickness in the small coupling regime. The situation in Proposition~\ref{p.thick} is more delicate since the curve of initial conditions (that is the meaning of $\Gamma_n(t)$) does not belong to the same level surface $S_I$ as in the discrete case. Nevertheless, the condition (3) from Proposition~\ref{p.thick} ensures that the value of the Fricke-Vogt invariant $I$ does not change too fast along the curve $\Gamma_n(t)$, and that in turn allows one to modify the main steps of the proof of \cite[Theorem~1.2]{DG2011} in order to derive Proposition~\ref{p.thick}. We carry out the details of the proof of Proposition~\ref{p.thick} in Subsections~\ref{s.s04}, \ref{ss.egs}, \ref{ss.prelimEst}, and \ref{ss.proofdist}.

\subsection{Choice of a Coordinate System in a Neighborhood of a Singular Point}\label{s.s04}

Due to the smoothness of the invariant manifolds of the curve of periodic points of period two described above, there exists a smooth change of coordinates $\Phi:O_{r_0}(P_1)\to \mathbb{R}^3$ such that $\Phi(P_1)=(0,0,0)$ and
\begin{itemize}

\item $\Phi(\mathfrak{I})$ is a part of the line $\{x=0, z=0\}$;

\item $\Phi(W^{\cs}_{\loc}(\mathfrak{I}))$ is a part of the plane $\{z=0\}$;

\item $\Phi(W^{\cu}_{\loc}(\mathfrak{I}))$ is a part of the plane $\{x=0\}$;

\item $\Phi(W^{\sss}_{\loc}(P_1))$ is a part of the line $\{y=0, z=0\}$;

\item $\Phi(W^{\uu}_{\loc}(P_1))$ is a part of the line $\{x=0, y=0\}$.

\end{itemize}
Denote $f = \Phi \circ T \circ \Phi^{-1}$. Then,
$$
A \equiv Df(0,0,0)=D(\Phi\circ T\circ \Phi^{-1})(0,0,0)=\begin{pmatrix}
                                       \xi^{-1} & 0 & 0 \\
                                       0 & -1 & 0 \\
                                       0 & 0 & \xi \\
                                     \end{pmatrix},
$$
where $\xi$ is the largest eigenvalue of the differential $DT(P_1):T_{P_1}\mathbb{R}^3\to T_{P_1}\mathbb{R}^3$,
$$
DT(P_1)=\begin{pmatrix}
                                       2 & 2 & -1 \\
                                       1 & 0 & 0 \\
                                       0 & 1 & 0 \\
                                     \end{pmatrix}, \ \ \ \xi=\frac{3+\sqrt{5}}{2}=\mu^2.
$$
Let us denote $\frak{S}_I=\Phi(S_I)$. Then, away from $(0,0,0)$, the family $\{\frak{S}_I\}$ is a smooth family of surfaces, $\frak{S}_0$ is diffeomorphic to a cone, contains the lines $\{y=0, z=0\}$ and $\{x=0, y=0\}$, and at each non-zero point on
these lines, it has a quadratic tangency with a horizontal or vertical plane.

Due to the symmetries of the trace map, similar changes of coordinates exist in a neighborhood of each of the other singularities. Denote $O_{r_0} = O_{r_0}(P_1) \cup O_{r_0}(P_2) \cup O_{r_0}(P_3) \cup O_{r_0}(P_4)$.
%

The next statement shows that, roughly speaking, if a point stays in a neighborhood where normalizing coordinates are defined for $N$ iterates, then it must be $\xi^{-N}$-close to the center-stable manifold of the curve of fixed points.

\begin{prop}[Proposition 3.4 from \cite{DG2011}]\label{p.dist}
Given $C_1>0, C_2>0, \xi>1$,
there exist $\delta_0 = \delta_0(C_1, C_2,  \xi)$, $N_0 = N_0(C_1, C_2, \xi, \delta_0) \in \mathbb{N}$, and $C^{**} > C^* > 0$ such that for any $\delta \in (0, \delta_0)$, the following holds.

Let $f : \mathbb{R}^3 \to \mathbb{R}^3$ be a $C^2$-diffeomorphism such that
\begin{itemize}

\item[{\rm (i)}] $\|f\|_{C^2}\le C_1$;

\item[{\rm (ii)}] the planes $\{z=0\}$ and $\{x=0\}$ are invariant under iterates of $f$;

\item[{\rm (iii)}] every point of the line $\{z=0, x=0\}$ is a fixed point of $f$;

\item[{\rm (iv)}] at a point $Q\in \{z=0, x=0\}$ we have
$$
Df(Q)=\begin{pmatrix}
    \xi^{-1} & 0 & 0 \\
    0 & 1 & 0 \\
    0 & 0 & \xi \\
  \end{pmatrix};
$$

\item[{\rm (v)}] $\|Df(p) - A\| < \delta$ for every $p \in \mathbb{R}^3$, where
$$
A = Df(Q) = \begin{pmatrix}
    \xi^{-1} & 0 & 0 \\
    0 & 1 & 0 \\
    0 & 0 & \xi \\
  \end{pmatrix}.
$$
\end{itemize}
Introduce the following cone fields in $\mathbb{R}^3$:
\begin{align}
\label{e.coneformula1} K_p^{\rm u} & = \{\mathbf{v}\in T_p\mathbb{R}^3, \ \mathbf{v}=\mathbf{v}_{xy}+\mathbf{v}_z : |\mathbf{v}_z|\ge C_2 \sqrt{|z_p|}|\mathbf{v}_{xy}|\}, \\
\label{e.cuconeformula} K_p^{\cu} & = \{\mathbf{v} \in T_p\mathbb{R}^3, \ \mathbf{v}=\mathbf{v}_{x}+\mathbf{v}_{yz} : |\mathbf{v}_x|<0.01\xi^{-1} |\mathbf{v}_{yz}| \}, \\
\label{e.coneformulas} K_p^{\rm s} & = \{\mathbf{v}\in T_p\mathbb{R}^3, \ \mathbf{v}=\mathbf{v}_{x}+\mathbf{v}_{yz} : |\mathbf{v}_x|\ge
C_2 \sqrt{|x_p|}|\mathbf{v}_{yz}|\}, \\
\label{e.csconeformula} K_p^{\cs} & = \{\mathbf{v}\in T_p\mathbb{R}^3, \ \mathbf{v}=\mathbf{v}_{z}+\mathbf{v}_{xy} : |\mathbf{v}_z|<0.01\xi^{-1} |\mathbf{v}_{xy}| \}.
\end{align}

Suppose that for a finite orbit $p_0, p_1, p_3, \ldots, p_N$, we have
$$
(p_0)_x\ge 1, \ (p_1)_x<1, \ (p_N)_z\ge 1, \ (p_{N-1})_z<1,
$$
and there are curves $\gamma_0$ and $\gamma_N$ such that $\gamma_0$ connects $p_0$ with $W^{\sss}(Q)$ and is tangent to both cone fields $K^{\rm u}$ and $K^{\cu}$, and $\gamma_N$ connects $p_N$ with $W^{\uu}(Q)$ and is tangent to both
cone fields $K^{\rm s}$ and $K^{\cs}$.

Then
\begin{align*}
C^*\xi^{-N} & \le |(p_0)_z|\le C^{**}\xi^{-N}, \ \ \text{\rm
and} \\
C^*\xi^{-N} & \le |(p_N)_x|\le C^{**}\xi^{-N}.
\end{align*}
\end{prop}

From the proof of Proposition \ref{p.dist} (i.e., of \cite[Proposition 3.4]{DG2011}) one can extract the following statement:

\begin{lemma}\label{l.bk}
In the setting of Proposition \ref{p.dist}, denote $b_k=\text{\rm dist}(p_k, Q)$. There exists a constant $C'>0$ (independent of $k$ or $N$) such that the following holds:

If $k<N/2$, then $b_k\le C'(\xi-\delta)^{-k}$.

If $k\ge N/2$, then $b_k\le C'(\xi-\delta)^{-N+k}$.
\end{lemma}

Other statements from \cite{DG2011} we will need are the following:

\begin{prop}[Proposition 3.15 from \cite{DG2011}]\label{p.vectors} Given $C_1>0, C_2>0, \lambda>1$,
there exist $\delta_0=\delta_0(C_1, C_2,  \lambda)$, $N_0\in \mathbb{N}, N_0=N_0(C_1, C_2,  \lambda,  \delta_0)\in
\mathbb{N}$, and $\widetilde{C}>0$ such that for any $\delta\in (0, \delta_0)$, the
following holds.

Under the conditions of and with the notation from Proposition~\ref{p.dist}, suppose that $\mathbf{v} \in T_{p_0}\mathbb{R}^3, \mathbf{v}\in K^u_{p_0}$. Then $|Df_{p_0}^N(\mathbf{v})| \ge \widetilde{C}\xi^{N/2}|\mathbf{v}|$.
\end{prop}

Let us denote $\mathbf{v}_k=Df^k(\mathbf{v}), k=0, 1, \ldots, N,$ and $D_k=|(\mathbf{v}_k)_z|, d_k=|(\mathbf{v}_k)_{xy}|$.

\begin{lemma}[Lemma 3.16 from \cite{DG2011}]\label{l.kstar}
There exists $k^*$ such that $d_k\ge D_k$ for all $k\le k^*$, and $d_k<D_k$ for all $k>k^*$.
\end{lemma}

\begin{lemma}[Lemma 3.17 from \cite{DG2011}]\label{l.onkstarboundabove}
There is a constant $C_9$ independent of $N$ such that
$$
\xi^{k^*}\le C_9\xi^{N/2}.
$$
\end{lemma}

\begin{lemma}\label{l.c6}
There exists $C_6>0$ such that for all $k\le k^*$ we have
$$
C_6^{-1}\le |(\mathbf{v}_k)_y|\le C_6.
$$
\end{lemma}

\begin{proof}[Proof of Lemma \ref{l.c6}]
From the proof of Lemma \ref{l.kstar} (i.e., \cite[Lemma 3.16]{DG2011}) one can immediately extract that $|(\mathbf{v}_k)_y|\le d_k\le C_6$ for some constant $C_6>0$ for all $k\le k^*$. At the same time, since $\|f\|_{C^2}\le C_1$, and
$$
A = Df(Q) = \begin{pmatrix}
    \xi^{-1} & 0 & 0 \\
    0 & 1 & 0 \\
    0 & 0 & \xi \\
  \end{pmatrix},
$$
we have
$$
|(\mathbf{v}_{k+1})_y|\ge |(\mathbf{v}_k)_y|-\min{(\delta, Cb_k)}\|\mathbf{v}_k\|\ge (1-C^*b_k)|(\mathbf{v}_k)_y|,
$$
where $C^*$ is a constant independent of $k<k^*$ and $N$. Hence, taking into account Lemma \ref{l.bk}, we have
$$
|(\mathbf{v}_k)_y|\ge \left[\prod_{i=1}^k(1-C^*b_k)\right]|(\mathbf{v}_0)_y|\ge C^{-1}_6
$$
 for some constant $C_6>0$.
\end{proof}

\begin{prop}[Proposition 3.18 from \cite{DG2011}]\label{p.distances}
Given $C_1 > 0, C_2 > 0, C_3 > 0, \xi > 1$,
there exist $\delta_0 = \delta_0(C_1, C_2, C_3, \xi)$, $N_0 = N_0(C_1, C_2, C_3, \lambda,  \delta_0) \in \mathbb{N}$, and $C > 0$ such that for any $\delta \in (0, \delta_0)$ and any $\Delta > 0$, 
the following holds.

Under the conditions of and with the notation from Proposition~\ref{p.dist}, suppose that the curve $\gamma_0$ has a
curvature bounded by $C_3$. Suppose also that for the points $p = (x_p, y_p, z_p)$ and $q = (x_q, y_q, z_q)$, the following holds:

\begin{enumerate}

\item $p, q\in \gamma_0$;

\item For some $N\ge N_0$ both $f^N(p)$ and $f^N(q)$ have $z$-coordinates larger than 1, and both $f^{N-1}(p)$ and $f^{N-1}(q)$ have $z$-coordinates not greater than 1;

\item $\mathrm{dist}(f^N(p), f^N(q)) = \Delta$.

\end{enumerate}

Denote $p_k=f^k(p)$, $q_k=f^k(q)$, $k=0, \ldots, N$. Let $\mathbf{v}\in T_p\mathbb{R}^3$ and  $\mathbf{w}\in
T_q\mathbb{R}^3$ be vectors tangent to $\gamma_0$ and denote $\mathbf{v}_k = Df^k(v)$ and $\mathbf{w}_k = Df^k(w)$, $k = 0, \ldots, N$. Let $\alpha_k$ be the angle between $\mathbf{v}_k$ and $\mathbf{w}_k$.

Then,
\begin{equation}\label{e.sumsarebounded}
\sum_{k=0}^{N} \alpha_k < C \Delta \; \text{ and } \; \sum_{k=0}^{N} \mathrm{dist}(p_k, q_k) < C\Delta.
\end{equation}
\end{prop}

We will need to modify Proposition \ref{p.distances}.

\begin{defi}
For any points $p, q$ and any nonzero vectors $\mathbf{v}_p, \mathbf{v}_q$, define
\begin{equation}\label{e.ffunction}
\mathfrak{F}(p, q, \mathbf{v}_p, \mathbf{v}_q) \equiv \frac{\angle (\mathbf{v}_p, \mathbf{v}_q)}{\text{\rm dist}(p,q)}.
\end{equation}
\end{defi}

\begin{defi}\label{d.fgamma}
For any $C^2$ smooth curve $\gamma$, define
\begin{equation}\label{e.ffunction}
\mathfrak{F}(\gamma) = \max \mathfrak{F}(p, q, \mathbf{v}_p, \mathbf{v}_q),
\end{equation}
where  $\mathbf{v}_p$ and $\mathbf{v}_q$ are unit vectors tangent to $\gamma$ at the points $p$ and $q$, and the maximum is taken over all pairs of distinct points $p, q\in \gamma$.
\end{defi}

\begin{remark}\label{r.curv}
Notice that $\mathfrak{F}(\gamma)$ is the maximum of the curvature of the curve $\gamma$ over all its points.
\end{remark}

Here is the formal statement that we need:

\begin{prop}[Modified version of Proposition \ref{p.distances}]\label{p.distancesmod}
Given $C_1 > 0, C_2 > 0,  \xi > 1$,
there exist $\delta_0 = \delta_0(C_1, C_2,  \xi)$, $N_0 = N_0(C_1, C_2,  \xi,  \delta_0) \in \mathbb{N}$, and $C > 0$ such that for any $\delta \in (0, \delta_0)$ and any $\Delta > 0$, 
the following holds.

Under the conditions of and with the notation from Proposition~\ref{p.dist}, suppose that for the points $p = (x_p, y_p, z_p)$ and $q = (x_q, y_q, z_q)$, the following holds:

\begin{enumerate}

\item $p, q\in \gamma_0$;

\item for some $N \ge N_0$, both $f^N(p)$ and $f^N(q)$ have $z$-coordinates larger than 1, and both $f^{N-1}(p)$ and $f^{N-1}(q)$ have $z$-coordinates not greater than 1;

\item $\mathrm{dist}(f^N(p), f^N(q)) = \Delta$.

\end{enumerate}

Denote $p_k = f^k(p)$, $q_k = f^k(q)$, $k = 0, \ldots, N$. Let $\mathbf{v} \in T_p \mathbb{R}^3$ and  $\mathbf{w} \in T_q \mathbb{R}^3$ be vectors tangent to $\gamma_0$ and denote $\mathbf{v}_k = Df^k(\mathbf{v})$ and $\mathbf{w}_k = Df^k(\mathbf{w})$, $k = 0, \ldots, N$. Let $\alpha_k$ be the angle between $\mathbf{v}_k$ and $\mathbf{w}_k$.

Denote $\gamma_k = f^k(\gamma_0)$, and set $\mathfrak{F}_k = \mathfrak{F}(\gamma_k)$, $k = 0, \ldots, N$.

Then,
\begin{equation}\label{e.sumsarebounded}
\sum_{k=0}^{N} \mathrm{dist}(p_k, q_k) < C\Delta,   \ \ \text{ and } 
\end{equation}
\begin{equation}\label{e.curvature}
 \mathfrak{F}_N< \xi^{-N/4}\mathfrak{F}_0+C.   
\end{equation}

Moreover, for any $C_3>0$ there exists $\hat C>0$ such that if $\mathfrak{F}_0<C_3$, then
\begin{equation}\label{e.sumofanglessarebounded}
\sum_{k=0}^{N} \alpha_k < \hat C \Delta.
\end{equation}
\end{prop}

\begin{proof}[Proof of Proposition  \ref{p.distancesmod}]
The inequalities (\ref{e.sumsarebounded}) and (\ref{e.sumofanglessarebounded}) follow from Proposition \ref{p.distances}, one just needs to observe that in the proof of Proposition \ref{p.distances} (i.e., of \cite[Proposition 3.18]{DG2011})  the estimates on $\sum_{k=0}^{N} \mathrm{dist}(p_k, q_k)$ do not use any bound on the curvature of $\gamma_0$. Hence we only need to prove (\ref{e.curvature}).

Notice that it is enough to prove Proposition~\ref{p.distancesmod} in the case when the points $p$ and $q$ are arbitrarily close to each other. This follows from the fact that in Definition~\ref{d.fgamma} one can take the maximum over those pairs of points that are sufficiently close, compare with Remark~\ref{r.curv}.

Denote by $\Gamma$ the piece of the curve $\gamma_0$ between the
points $p_0$ and $q_0$, and set $\Gamma_k=f^k(\Gamma)$, $k=0, 1,
2, \ldots, N$. Denote $\mu_k=|\Gamma_k|$. Due to the remark above
we can assume that $\mu_0\ll 1$, and that for any vector tangent to $\Gamma$, the value
of $k^*$ is the same.

We will need to use the following elementary estimates, whose proofs we leave to the reader.

\begin{lemma}\label{l.basic1}
Given $\alpha_0>0$, $\tau>1$, $\delta>0$, suppose that
$$
\alpha_{k+1}\le \tau\alpha_{k}+\delta, \ \ k=0, 1, \ldots
$$
Then $\alpha_k\le \tau^k\left(\alpha_0+\frac{\delta}{\tau-1}\right)$.
\end{lemma}

\begin{lemma}\label{l.basic2}
Given $\alpha_0>0$, $0<t<1$, $\Delta>0$, suppose that
$$
\alpha_{k+1}\le t\alpha_{k}+\Delta, \ \ k=0, 1, \ldots
$$
Then $\alpha_k\le t^k\alpha_0+\frac{\Delta}{1-t}$.
\end{lemma}

Recall that $\mathbf{v}_k$ and $\mathbf{w}_k$ are the vectors tangent to $\Gamma_k$ at $p_k=f^k(p_0)$ and $q_k=f^k(q_0)$, and $\alpha_k=\angle(\mathbf{v}_k, \mathbf{w}_k)$.

We know that for $k \le k^*$, one has
$$
\alpha_{k+1} \le (\xi + \delta) \alpha_k + C_1 \mu_k \le (\xi + \delta) \alpha_k + C_1 C_6 \mu_0,
$$
since Lemma~\ref{l.c6} is applicable in this case.

We also know that for some $m \in \mathbb{N}$ independent of $N$, and for all $k > k^* + m$, we have
$$
\alpha_{k+1} \le (\xi^{-1} + \delta) \alpha_k + C_1 \mu_k \le (\xi^{-1} + \delta) \alpha_k + C_1\mu_N, \ \ \text{\rm and}\ \ \ \mu_{k+1}\ge (\xi-2\delta)\mu_k.
$$
Therefore, Lemma~\ref{l.basic1} gives
$$
\alpha_{k^*} \le (\xi + \delta)^{k^*} \left( \alpha_0 + \frac{C_1 C_6 \mu_0}{\xi + \delta-1} \right).
$$
We also have
$$
\alpha_{k^* + m} \le C_1^m \alpha_{k^*} + m C_1^m \mu_0.
$$
Therefore, due to Lemma~\ref{l.basic2} we have
$$
\alpha_N \le (\xi^{-1} + 2 \delta)^{N/2} \left( C_1^m \alpha_{k^*} + m C_1^m \mu_0 \right) + \frac{\mu_N}{1 - \xi^{-1} - 2 \delta}.
$$
Combining these inequalities, we get
\begin{multline*}
    \frak{F}(p_N, q_N, \mathbf{v}_N, \mathbf{w}_N) = \frac{\alpha_N}{\mu_N} \le \\
    \const \cdot \frac{(\xi^{-1} + 2 \delta)^{N/2}(\xi + \delta)^{k^*}}{(\xi - 2 \delta)^{N/2}} \frak{F}(p_0, q_0, \mathbf{v}_0, \mathbf{w}_0) + \const \cdot \frac{(\xi^{-1} + 2 \delta)^{N/2}}{(\xi - 2 \delta)^{N/2}} + \const,
\end{multline*}
and \eqref{e.curvature} follows.
\end{proof}


\subsection{Distortion Property: Estimate of the Gap Sizes}\label{ss.egs}

Let us consider $\Gamma_n$. Its image $T^{6m}(\Gamma_n)$ intersects $R$ in a finite number of closed curves. The gaps between these curves correspond to the gaps in the parameter space $J_n$ that we will call the gaps of order one. Denote
$$
M^-_n = \min_{t \in J_n} I(\Gamma_n(t)), \ \ M^+_n = \max_{t \in J_n} I(\Gamma_n(t)).
$$
The assumption $(1)$ from Proposition~\ref{p.thick} implies that $M^-_n, M^+_n\to 0$ as $n\to \infty$, and assumption $(3)$ implies that for some uniform constant $C_1>1$ that is independent of $n$, the size of any gap of order one is between $C_1^{-1} \sqrt{M^-_n}$ and $C_1 \sqrt{M^+_n}$.

Those gaps in $J_n$ that are formed by the intersection of $T^{6m+(k-1)}(\Gamma_n)$ with the complement of $R$ and do not have order less than $k$ will be called gaps of order $k$. It is clear that every gap in $\Gamma_n \cap
\Omega^+$ has some finite order. Therefore we have ordered all the gaps.

Consider some gap ${\mathfrak{l}}_G \subset J_n$ of order $k$. A \textit{bridge} that corresponds to this gap is a connected component of the complement of the union of all gaps of order $\le k$ next to the gap. There are two bridges that correspond to the chosen gap, take one of them, and denote it by $\mathfrak{l}_B$. Now let us consider $\mathfrak{L}_G \equiv \Gamma_n^{-1} \circ T^{6m+(k-1)}(\Gamma_n(\mathfrak{l}_G))$ and ${\mathfrak{L}}_B \equiv
\Gamma_n^{-1} \circ T^{6m+(k-1)}(\Gamma_n(\mathfrak{l}_B))$. By the definition of the order $k$ of the gap we know that
\begin{equation}\label{e.impeq}
C_3 \sqrt{M^-_n} \le \frac{|\mathfrak{L}_G|}{|\mathfrak{L}_B|} \le C_4 \sqrt{M^+_n}
\end{equation}
for suitable constants $C_3$ and $C_4$ independent of $n$.

\begin{prop}[Analog of Proposition~3.11 from \cite{DG2011}]\label{p.distortion}
There is a constant $K > 1$ independent of the choice of the gap and of $n$ such that
$$
K^{-1} \frac{|\mathfrak{l}_G|}{|\mathfrak{l}_B|} \le \frac{|\mathfrak{L}_G|}{|\mathfrak{L}_B|} \le K \frac{|\mathfrak{l}_G|}{|\mathfrak{l}_B|}.
$$
\end{prop}

After some preparatory work in Subsection~\ref{ss.prelimEst}, we will prove Proposition~\ref{p.distortion} in Subsection~\ref{ss.proofdist}. In the meantime, let us observe that the distortion property elucidated in Proposition~\ref{p.distortion} suffices to establish Proposition~\ref{p.thick}.

\begin{proof}[Proof of Proposition~\ref{p.thick}]
The statement follows from Proposition~\ref{p.distortion}, inequality (\ref{e.impeq}), and the fact that $M^-_n, M^+_n \to 0$ as $n \to \infty$.
\end{proof}

\subsection{Preliminary Estimates} \label{ss.prelimEst}

Let us denote by $U_\nu$ the bounded connected component of $\left( \bigcup_{0 \le I \le \nu} S_I \right) \backslash O_{r_2}$. Let us denote by $\pi_I$ the orthogonal projection of $S_I \cap U_\nu$ to $\mathbb{S} \backslash O_{r_2}$. Notice that $F^{-1}(\mathbb{S} \backslash O_{r_2})$ is a torus without small neighborhoods of the preimages of the singularities. We can define
$$
F_{[0, \nu]}:U_\nu\to \mathbb{T}^2\times [0, \nu]
$$
by
$$
F_{[0, \nu]}(x)=(F^{-1}(\pi_{I(x)}(x)), I(x)),
$$
and
$$
\widetilde T: F_{[0, \nu]}(U_\nu)\to \mathbb{T}^2\times [0, \nu]
$$
by
$$
\widetilde T(\widetilde x, I)=(F^{-1}\circ \pi_I\circ T\circ \pi_I^{-1}\circ F(\widetilde x), I).
$$
In this case $\widetilde T$ 
is $C^2$-close to $\mathcal{A} \times id : \mathbb{T}^2 \times [0,\mu] \to \mathbb{T}^2\times [0,\mu]$ if $\nu$ is small, and
$$
\widetilde T = H^{-1} \circ T \circ H,
$$
where $H(x, I) = \pi_{I(x)}^{-1} \circ F(x)$, that is, $\widetilde T$ and $T$ are smoothly conjugate. Set
\begin{equation}\label{e.ctilda}
\widetilde C := \max (\|H\|_{C^2}, \|H^{-1}\|_{C^2}, \|\Phi\|_{C^2}, \|\Phi^{-1}\|_{C^2}, \|H\circ\Phi\|_{C^2}, \|(H\circ\Phi)^{-1}\|_{C^2}).
\end{equation}
Informally speaking, the constant $\widetilde C$ gives an upper bound on the distortion induced by any change of coordinates that we may want to consider.

Define a cone field in $\mathbb{T}^2\times [0, \nu]$ using stable-unstable directions of the map $\mathcal{A}$ as $x$ and $y$ coordinates:
$$
K^{\cu}_{(x, y, I)}=\{\bar v=(v_x, v_y, v_I)\ |\ |v_y|>100|v_x|, \ |v_I|<C_0I\}.
$$
Due to assumption (3) in Proposition \ref{p.thick} and the fact that the curve $\Pi$ is transversal to the stable foliation on $\mathbb{S}$, one can choose $C_0>0$ and $k_0\in \mathbb{N}$ in such a way that $\widetilde T^k\circ H^{-1}(\Gamma_n)$ is tangent to $K^{\cu}$ for all $k\ge k_0$.

\begin{lemma}[Analog of Lemma 3.20 from \cite{DG2011}]\label{l.angels}
For $\nu > 0$ small enough, there exists $\eta \in (0,1)$ such that for any $p,q \in F_{[0, \nu]}(U_\nu)$ and unit vectors $\mathbf{v}_p\in K_p^{\cu}$, $\mathbf{v}_q\in K_q^{\cu}$, we have
$$
\angle (D \widetilde{T}_{p}(\mathbf{v}_p), D \widetilde{T}_{q}(\mathbf{v}_q)) \le \eta \angle(\mathbf{v}_p, \mathbf{v}_q) + 2 \|\widetilde{T}\|_{C^2} \text{\rm dist}(p,q).
$$
\end{lemma}

\begin{proof}[Proof of Lemma \ref{l.angels}]
If $\nu$ is small, then $\widetilde{T}$ is $C^2$-close to the map $\mathcal{A} \times id$. In particular, for any point $p \in F_{[0, \nu]}(U_\nu)$ and any vectors $\mathbf{v}_1, \mathbf{v}_2 \in K^{\cu}_p$,
$$
\angle (D\widetilde{T}_{p}(\mathbf{v}_1), D\widetilde{T}_{p}(\mathbf{v}_2)) \le \eta\angle(\mathbf{v}_1, \mathbf{v}_2),
$$
where $\eta \in (0,1)$ can be chosen uniformly for all  $p\in F_{[0, \nu]}(U_\nu)$. Therefore we
have
\begin{align*}
\angle (D\widetilde{T}_{p}(\mathbf{v}_p), D\widetilde{T}_{q}(\mathbf{v}_q)) & \le
\angle (D\widetilde{T}_{p}(\mathbf{v}_p), D\widetilde{T}_{p}(\mathbf{v}_q)) + \angle
(D\widetilde{T}_{p}(\mathbf{v}_q), D \widetilde{T}_{q}(\mathbf{v}_q)) \\
& \le \eta\angle(\mathbf{v}_p, \mathbf{v}_q)+2\|D\widetilde{T}_{p}(\mathbf{v}_q)-D\widetilde{T}_{q}(\mathbf{v}_q)\| \\
& \le \eta\angle(\mathbf{v}_p, \mathbf{v}_q)+2\|\widetilde{T}\|_{C^2}\text{\rm
dist}(p,q),
\end{align*}
as claimed.
\end{proof}

\begin{lemma}[Analog of Lemma 3.22 from \cite{DG2011}]\label{l.inequality} For $p,q\in  F_{[0, \nu]}(U_\nu)$, $p\ne q$, and vectors $\mathbf{v}_p\in K_p^{\cu}$, $\mathbf{v}_q\in
K_q^{\cu}$, consider the function $\mathfrak{F}(p, q, \mathbf{v}_p, \mathbf{v}_q)$. 
Suppose that $p$ and $q$ belong to a curve that is tangent to the cone field $K^{\cu}$. Then
$$
\frak{F}(\widetilde{T}(p),\widetilde{T}(q), D\widetilde{T}_{p}(\mathbf{v}_p), D\widetilde{T}_{q}(\mathbf{v}_q))\le \eta\mathfrak{F}(p, q, \mathbf{v}_p, \mathbf{v}_q) +2\|\widetilde{T}\|_{C^2}.
$$
In particular, if $\mathfrak{F}(p, q, \mathbf{v}_p, \mathbf{v}_q) >
\frac{4\|\widetilde{T}\|_{C^2}}{1-\eta}$, then
$$
\frak{F}(\widetilde{T}(p),\widetilde{T}(q), D\widetilde{T}_{p}(\mathbf{v}_p), D\widetilde{T}_{q}(\mathbf{v}_q))\le
\frac{1+\eta}{2}\mathfrak{F}(p, q, \mathbf{v}_p, \mathbf{v}_q).
$$
\end{lemma}
\begin{proof}[Proof of Lemma \ref{l.inequality}]
We have
\begin{align*}
\frak{F}(\widetilde{T}(p),\widetilde{T}(q), D\widetilde{T}_{p}(\mathbf{v}_p), D\widetilde{T}_{q}(\mathbf{v}_q)) & =\frac{\angle(D\widetilde{T}_{p}(\mathbf{v}_p),
D\widetilde{T}_{q}(\mathbf{v}_q))}{\text{\rm dist}(\widetilde{T}(p),\widetilde{T}(q))} \\
& \le \frac{\eta\angle(\mathbf{v}_p, \mathbf{v}_q)+2\|\widetilde{T}\|_{C^2}\text{\rm dist}(p,q)}{\text{\rm dist}(p,q)} \\
& = \eta\mathfrak{F}(p, q, \mathbf{v}_p, \mathbf{v}_q) +2\|\widetilde{T}\|_{C^2}.
\end{align*}
If we also have $\mathfrak{F}(p, q, \mathbf{v}_p, \mathbf{v}_q) > \frac{4\|\widetilde{T}\|_{C^2}}{1-\eta}$, then 
\begin{align*}
\eta\mathfrak{F}(p, q, \mathbf{v}_p, \mathbf{v}_q) + 2 \|\widetilde{T}\|_{C^2} & \le
\eta\mathfrak{F}(p, q, \mathbf{v}_p, \mathbf{v}_q) +\frac{1-\eta}{2}\mathfrak{F}(p, q, \mathbf{v}_p, \mathbf{v}_q) \\
& = \frac{1+\eta}{2}\mathfrak{F}(p, q, \mathbf{v}_p, \mathbf{v}_q).
\end{align*}
\end{proof}

Denote $\Psi=\Phi\circ H$. From (\ref{e.ctilda}) we have:

\begin{lemma}[Analog of Lemma 3.5 from \cite{DG2011}]\label{l.inequalities} Let $\widetilde{C}>0$ be given by (\ref{e.ctilda}).  There is $\nu>0$
such that the following holds. Suppose that $a, b\in
F_{[0, \nu]}(U_\nu)$, $\mathbf{v}_a\in T_a\mathbb{T}^2$,
$\mathbf{v}_b\in T_b\mathbb{T}^2$. Then the following inequalities hold :
\begin{align*}
\text{\rm dist}(H(a),H(b)) & \le \widetilde{C} \, \text{\rm
dist}(a,b), \\
\text{\rm dist}(a,b) & \le \widetilde{C} \, \text{\rm
dist}(H(a),H(b)), \\
\angle(DH(\mathbf{v}_a), DH(\mathbf{v}_b)) & \le \widetilde{C} (\angle(\mathbf{v}_a,
\mathbf{v}_b)+\text{\rm dist}(a,b)), \\
\angle(\mathbf{v}_a, \mathbf{v}_b) & \le \widetilde{C} (\angle(DH(\mathbf{v}_a),
DH(\mathbf{v}_b))+\text{\rm dist}(H(a),H(b))).
\end{align*}
Moreover, if $\Psi(a)$ and $\Psi(b)$ are defined, then
\begin{align*}
\text{\rm dist}(\Psi(a),\Psi(b)) & \le \widetilde{C} \,
\text{\rm
dist}(a,b), \\
\text{\rm dist}(a,b) & \le \widetilde{C} \, \text{\rm
dist}(\Psi(a),\Psi(b)), \\
\angle(D\Psi(\mathbf{v}_a), D\Psi(\mathbf{v}_b)) & \le \widetilde{C}(\angle(\mathbf{v}_a,
\mathbf{v}_b)+\text{\rm dist}(a,b)), \\
\angle(\mathbf{v}_a, \mathbf{v}_b) & \le \widetilde{C}(\angle(D\Psi(\mathbf{v}_a),
D\Psi(\mathbf{v}_b))+\text{\rm dist}(\Psi(a),\Psi(b))).
\end{align*}
\end{lemma}

From Lemma \ref{l.inequalities} and the definition of $\frak{F}$ we get the following:

\begin{lemma}[Analog of Lemma 3.23 from \cite{DG2011}]\label{l.changeofcoordinates} Fix a small $\nu\ge 0$. Suppose
that $a, b \in F_{[0, \nu]}(U_\nu)$ are such
that $\Psi(a)$ and $\Psi(b)$ are defined, and $\mathbf{v}_a\in
T_a\mathbb{T}^2$, $\mathbf{v}_b\in
T_b\mathbb{T}^2$. Then,
$$
\frak{F}(a, b, \mathbf{v}_a, \mathbf{v}_b)\le \widetilde{C}^2(\frak{F}(\Psi(a),
\Psi(b), D\Psi(\mathbf{v}_a), D\Psi(\mathbf{v}_b))+1)
$$
and
$$
\frak{F}(\Psi(a), \Psi(b), D\Psi(\mathbf{v}_a), D\Psi(\mathbf{v}_b))\le \widetilde{C}^2(\frak{F}(a, b, \mathbf{v}_a, \mathbf{v}_b)+1).
$$
\end{lemma}

Now we are ready to choose a uniform upper bound on the curvatures of all the images of all the curves $\Gamma_n$ that are outside of small neighborhoods of the singularities. Namely, recall that the constant $C$ was given by \eqref{e.curvature}, $\widetilde C$ was defined by \eqref{e.ctilda}, and define
$$
M_{\frak{F}} := \max \left(2{\widetilde C}^2(1+C), \frac{40\|\widetilde T\|_{C^2}}{1-\eta}\right).
$$
Choose $N_0\in \mathbb{N}$ large enough to make sure that
$$
\xi^{-N/4}<\frac{1}{4{\widetilde C}^2}.
$$
Then for all large $n$ and for any small piece $\gamma$ of the curve $\Gamma_n$, all the iterates of the curve $\gamma$ that are in $U_\nu$ will have curvature not greater than $M_{\frak{F}}$. Indeed, after some finite number of initial iterates $k$, the curve $H^{-1}(T^k(\gamma))$ is tangent to the cone field $K^{\cu}$, will have curvature bounded by $M_{\frak{F}}$, and due to Lemma \ref{l.inequality} that curvature will remain bounded by $M_{\frak{F}}$ until we need to change the coordinates by applying the map $\Psi$. We will do that only if the image of $\gamma$ is going to spend more than $N_0$ iterates in a neighborhood of the singularity. After $N>N_0$ iterates in a neighborhood of the singularity, due to \eqref{e.curvature} the curvature will be bounded by
$$
\xi^{N/4}({\widetilde C}^2(M_{\frak{F}}+1))+C\le \frac{M_{\frak{F}}+1}{4{\widetilde C}^2}+C,
$$
and after application of the change of coordinates $\Psi^{-1}$ we get a curve with curvature bounded by
$$
\frac{1}{4}(M_{\frak{F}}+1) +{\widetilde C}^2C<M_{\frak{F}}.
$$

Let us now notice that the partial hyperbolicity of the map $\widetilde T$ together with \eqref{e.sumsarebounded} implies the following statement (the proof is similar to the proof of \cite[Lemma 3.24]{DG2011}, so we do not repeat it here):

\begin{lemma}[Analog of Lemma 3.24 from \cite{DG2011}]\label{l.vector}
There is $R_1>0$ such that the following holds for all sufficiently large $n$. Suppose that $\mathbf{v}$ is a non-zero vector tangent to $\Gamma_n$ at some point $p\in \Gamma_n$. Let $N\in \mathbb{N}$ be such that $T^N(p)$ belongs to the bounded component of $S_I\backslash O_{r_1}$, where $I=I(p)$. Then,
$$
\sum_{i=0}^{N}\|DT^i(\mathbf{v})\|\le R_1\|DT^N(\mathbf{v})\|.
$$
\end{lemma}

Lemma \ref{l.vector} implies the following statement.

\begin{lemma}\label{l.sum}
There are constants $R_1>0$ and $\kappa_1>0$ such that for all large $n$ and any  $N\in \mathbb{N}$, the following holds. Suppose that $\gamma\subset T^N(\Gamma_n)\backslash O_{r_1}$ is a connected curve of length not greater than $\kappa_1$. Let the points $p,q\in \Gamma_n$ be such that $T^N(p)\in \gamma$ and $T^N(q)\in \gamma$. Then,
$$
\sum_{i=0}^{N}\text{\rm dist}(T^i(p), T^i(q))<R_1.
$$
\end{lemma}

Finally, since we established a uniform bound on the curvature of all the images of the curves $\Gamma_n$ outside of small neighborhoods of singularities, we can use Lemma \ref{l.sum} together with \eqref{e.sumofanglessarebounded} to get the key technical statement that we need to establish the distortion property (i.e., Proposition~\ref{p.distortion}):

\begin{lemma}[Analog of Lemma 3.19 from \cite{DG2011}]\label{l.prelim}
There are constants $R>0$, and $\kappa>0$ such that for any large $n$ and $N\in \mathbb{N}$, the following holds. Suppose that $\gamma\subset T^N(\Gamma_n)\backslash O_{r_1}$ is a connected curve of length not greater than $\kappa$. Let the points $p,q\in \Gamma_n$ be such that $T^N(p)\in \gamma$ and $T^N(q)\in \gamma$, and $\mathbf{v}_p$ and $\mathbf{v}_q$ be unit vectors tangent to $\gamma$ at points $p$ and $q$. Then
$$
\sum_{i=0}^{N}\left(\angle(DT^i(\mathbf{v}_p), DT^i(\mathbf{v}_q))+\text{\rm dist}(T^i(p), T^i(q))\right)<R.
$$
\end{lemma}

\subsection{Proof of the Distortion Property}\label{ss.proofdist}

\begin{proof}[Proof of Proposition \ref{p.distortion}]

Notice that we need to prove that
$$
\left|\log \left(\frac{|\mathfrak{L}_G||\mathfrak{l}_B|}{|\mathfrak{L}_B||\mathfrak{l}_G|}\right)\right|
$$
is bounded by some constant independent of the choice of the index $n$ (as soon as $n$ is large enough) and the gap in $\Gamma_n$. There are points $p_G\in \mathfrak{l}_G$ and $p_B\in
\mathfrak{l}_B$ such that if $\mathbf{v}_G$ is a unit  vector tangent to the
curve $\mathfrak{l}_G$ at $p_G$, and $\mathbf{v}_B$ is a unit vector tangent to
the curve $\mathfrak{l}_B$ at $p_B$, then
\begin{align*}
\left|\log \left(\frac{|\mathfrak{L}_G||\mathfrak{l}_B|}{|\mathfrak{L}_B||\mathfrak{l}_G|}\right)\right| & =
\left|\log \left(\frac{|T^{n+2}(\mathfrak{l}_G)||\mathfrak{l}_B|}{|T^{n+2}(\mathfrak{l}_B)||\mathfrak{l}_G|}\right)\right| \\
& = \left|\log \left(\frac{|DT^{n+2}(\mathbf{v}_G)|}{|DT^{n+2}(\mathbf{v}_B)|}\right)\right| \\
& = \left|\sum_{i=0}^{n+1}\left(\log |DT|_{DT^i(\mathbf{v}_G)}(T^i(p_G))| - \log |DT|_{DT^i(\mathbf{v}_B)}(T^i(p_B))|\right)\right| \\
& \le \sum_{i=0}^{n+1}\left|\log |DT|_{DT^i(\mathbf{v}_G)}(T^i(p_G))| - \log |DT|_{DT^i(\mathbf{v}_B)}(T^i(p_B))|\right| \\
& \le \sum_{i=0}^{n+1}\left||DT|_{DT^i(\mathbf{v}_G)}(T^i(p_G))| -  |DT|_{DT^i(\mathbf{v}_B)}(T^i(p_B))|\right|.
\end{align*}

We estimate each of the terms in this sum using the following simple statement applied to the trace map $T$:

\begin{lemma}[Lemma 3.26 from \cite{DG2011}]\label{l.lastest}
Suppose $f : \mathbb{R}^n \to \mathbb{R}^n$ is a smooth map, $a, b\in \mathbb{R}^n$, and $\mathbf{v}_a \in T_a \mathbb{R}^n, \mathbf{v}_b \in T_b\mathbb{R}^n$ are unit vectors. Then,
$$
\left||Df|_{\mathbf{v}_a}(a)|-|Df|_{\mathbf{v}_b}(b)|\right|\le \|f\|_{C^2}(\angle(\mathbf{v}_a,\mathbf{v}_b)+|a-b|).
$$
\end{lemma}

Now Proposition \ref{p.distortion} follows from Lemma~\ref{l.prelim}.
\end{proof}

\subsection{Verifying the Assumptions of Propositions~\ref{p.thick}}

Here, we show how to verify the assumptions of Proposition~\ref{p.thick} in the models under consideration.

\begin{proof}[Proof of Theorem~\ref{t:hi}]
Denote the curve of initial conditions by
\[
\Gamma(t) = (u_1(t),u_0(t),u_{-1}(t)) \eqdef (x_1(t^2),x_0(t^2),x_{-1}(t^2)).
\]
We then choose $J_n = [2\pi n + \alpha_n, (2n+1)\pi - \beta_n]$ as above and denote by $\Gamma_n = \Gamma|_{J_n}$. We proceed by verifying that conditions (1)--(4) of Proposition~\ref{p.thick} hold true.
\medskip

\noindent \textbf{Property (1).} Write  $v_{-1}(t)=v_0(t) = \cos(t)$ and $v_1(t) = \cos(2t)$. Then, put $u_{j,n} = u_j|_{J_n}, v_{j,n} = v_j|_{J_n}$. First, by \cite[Theorem~1.3]{PoTru1987},
\[
\|u_{j,n} - v_{j,n}\|_\infty
\lesssim
n^{-1},
\]
showing that $u_{j,n}- v_{j,n} \to 0$ uniformly. In fact, the estimates from \cite{PoTru1987} hold for complex $t$ as well. In particular, since all functions in sight are analytic, $u_{j,n}^{(k)} - v_{j,n}^{(k)} \to 0$ uniformly for each $k$ (in particular for $k=1,2$) by Cauchy estimates.
\medskip

\noindent \textbf{Property (2).} This is immediate from the choice of $J_n$ and \eqref{e.1}.
\medskip

\noindent \textbf{Property (3).}  Note that
\[
\frac{d}{dt} \log I(\Gamma_n(t)) = 2\cot t- 2t^{-1} +  \frac{2 t \cot\sqrt{t^2-\lambda}}{\sqrt{t^2-\lambda}} - \frac{2t}{\sqrt{t^2-\lambda}}
\]
By our choice of $J_n$, this is uniformly bounded for large $n$.
\medskip

\noindent \textbf{Property (4).} This follows immediately from the discussion preceding Proposition~\ref{p.thick}.
\medskip

Consequently, the result follows from Proposition~\ref{p.thick} and Lemma~\ref{l.ABSC}.
\end{proof}

\begin{remark}\label{r.difficultgencase0}
In the proof of Theorem~\ref{t:hi} we made crucial use of the explicit expressions for the traces and the invariant in the case of constant $f_\a, f_\b$. It would of course be of interest to verify the assumptions of Proposition~\ref{p.thick} for more general choices of $f_\a, f_\b$. For example, one could in this way produce continuous or even smooth potentials. Alas, such an extension will not at all be straightforward. A modified version of property (1) and property (4) can be established in the general case, in fact by the exact same argument, 
and we expect that property (2) can be established under weaker assumptions as well. However, our proof of property (3)  relies on the explicit formula we have for the invariant in the locally constant case, and checking property (3) in a more general case seems to be a challenging technical task. 
\end{remark}

\section{The Low-Energy Regime} \label{s.lo}

Our goal in this section is to prove Theorem~\ref{t:lo}, which states that in the low-energy regime, the spectrum of the separable 2D Schr\"odinger operators we consider is a Cantor set. In this proof the existing work on the Fibonacci Hamiltonian will again be helpful, and indeed almost sufficient. As a consequence this section will be less demanding than Section~\ref{s.hi}.

Our goal is to show that by fixing a suitable compact energy interval $[0, E_0]$, $E_0 > 0$, and taking suitable building blocks, we can ensure that
\begin{equation}\label{e.smengoal1}
[0, E_0/2] \cap \Sigma \not= \emptyset
\end{equation}
and
\begin{equation}\label{e.smengoal2}
[0, E_0] \cap \Sigma = \bigcup_{k=0}^\infty A_k, \ \ \text{\rm where $A_k$ is compact, with}\ \  \dim_\mathrm{B}^+A_k < \frac12.
\end{equation}
It is then not hard to verify that
\begin{equation}\label{e.smencons1}
[0, E_0] \cap (\Sigma + \Sigma) \not= \emptyset
\end{equation}
and
\begin{equation}\label{e.smencons2}
\Leb([0, E_0] \cap (\Sigma + \Sigma)) = 0.
\end{equation}
Since we know that $\Sigma$ is a Cantor set, it then follows that $[0, E_0] \cap (\Sigma + \Sigma)$ is compact and nowhere dense. It could happen that $E_0$ is the edge of a gap of $\Sigma+\Sigma$ and hence is an isolated point of $[0,E_0] \cap (\Sigma+\Sigma)$. However, since $\Sigma$ itself has no isolated points, this is the only possible isolated point of $[0,E_0]\cap(\Sigma+\Sigma)$, and it is easy to avoid this by perturbing $E_0$. We will also ensure that $(-\infty, 0) \cap \Sigma = \emptyset$, which gives $(-\infty, 0) \cap (\Sigma + \Sigma) = \emptyset$. Thus, we will have accomplished our goal in the low-energy region.

Since the primary aim of this paper is to find examples that exhibit certain spectral phenomena, we choose our building blocks in such a way that the proofs of \eqref{e.smengoal1} and \eqref{e.smengoal2} are as simple as possible. We expect that the same properties can be shown in greater generality, but only after overcoming significant technical obstacles. To keep the length of the paper in check, we will not pursue this in detail here, but merely point out some of these obstacles in Remark~\ref{r.difficultgencase} below.

One of our building blocks will simply be the zero potential on some finite interval. This allows us to use test functions that are supported in this interval in order to produce low energies in the spectrum, and in particular verify \eqref{e.smengoal1} if $E_0$ and the interval length are compatible. Ensuring this is easy: as soon as one of them is fixed, the other one can be determined accordingly.

The choice of the second building block determines the difficulty of proving \eqref{e.smengoal2}. One needs to choose it so that the invariant restricted to $[0,E_0] \cap \Sigma$ is large. Recall that if the second building block consists of a constant potential, then the invariant can be computed explicitly. Moreover, in this particular case \cite{DFG2014} had already proved the desired statement, provided the constant value is sufficiently large. We will recall this derivation below. If, on the other hand, the second building block does not consist of a constant potential, then the invariant can in general not be determined explicitly, and the arguments used in \cite{DFG2014} to show that it is large on $[0,E_0] \cap \Sigma$ (which rely on the explicit formula) cannot be mimicked. In fact, it is not even clear in general how to proceed; see again Remark~\ref{r.difficultgencase} for further comments on the difficulties in establishing \eqref{e.smengoal2} in the general case.

Let us begin with \eqref{e.smengoal1}.

\begin{lemma} \label{l.smengoal1}
Suppose $\ell_\a > 0$ and $f_\a \in L^2(0,\ell_\a)$ is the zero element. Set
\begin{equation}\label{e.e0def}
E_0 := \frac{24}{\ell_\a^2}.
\end{equation}
Then \eqref{e.smengoal1} holds for any choice of $\ell_\b > 0$ and any non-negative $f_\b \in L^2(0,\ell_\b)$.
\end{lemma}

\begin{proof}
Let $E_0$ be given by \eqref{e.e0def} and suppose that $\ell_\b > 0$ and $f_\b \in L^2[0,\ell_\b)$ is non-negative. Fix any $\omega \in \Omega$. By \eqref{e.sigmaisconstant}, we have $\sigma(H_\omega) = \Sigma$. Thus, our goal is to show that
\begin{equation}\label{e.minmaxlemmagoal}
[0, E_0/2] \cap \sigma(H_\omega) \not= \emptyset.
\end{equation}
We consider the case where $\omega_0 = \a$; the other case is handled similarly (by considering the interval $[\ell_\b, \ell_\b + \ell_\a)$ instead of the interval $[0,\ell_\a)$).

Since $\omega_0 = \a$, $V_\omega$ vanishes on $[0,\ell_\a)$. Moreover, the function
$$
\varphi(x) = \begin{cases} \frac{2x}{\ell_\a} & x \in [0,\frac{\ell_\a}{2}) \\ 2 - \frac{2x}{\ell_\a} & x \in [\frac{\ell_\a}{2}, \ell_\a) \\ 0 & \text{otherwise} \end{cases}
$$
belongs to the form domain of $H_\omega$. We have
$$
\langle \varphi, H_\omega \varphi \rangle = \int_0^{\ell_\a} |\varphi'(x)|^2 \, dx = \ell_\a \cdot \frac{4}{\ell_\a^2} = \frac{4}{\ell_\a}.
$$
and
$$
\|\varphi\|^2 = \int_0^{\ell_\a} |\varphi(x)|^2 \, dx = 2 \int_0^{\frac{\ell_\a}{2}} \frac{4x^2}{\ell_\a^2} \, dx = \frac{8}{3 \ell_\a^2} \left( \frac{\ell_\a}{2} \right)^3 = \frac{\ell_\a}{3}.
$$
Thus,
$$
\frac{\langle \varphi, H_\omega \varphi \rangle}{\|\varphi\|^2} = \frac{12}{\ell_\a^2},
$$
and by the min-max principle (in the form version), it follows that $\sigma(H_\omega) \cap (-\infty, \frac{12}{\ell_\a^2}] \not= \emptyset$. On the other hand, since $H_\omega \ge 0$, we also know that $\sigma(H_\omega) \cap (-\infty, 0) = \emptyset$. Together with \eqref{e.e0def} this implies \eqref{e.minmaxlemmagoal}.
\end{proof}

Let us now turn our attention to \eqref{e.smengoal2}. As was mentioned above, we do not strive for maximal generality, but rather consider a case where the desired statement can be proved without the need to overcome serious technical difficulties.

\begin{lemma} \label{l.smengoal2}
Let $\ell_\a = \ell_\b = 1$, $f_\a  = 0 \cdot \chi_{(0,1)}$, and $f_\b = \lambda \cdot \chi_{(0,1)}$. Then, \eqref{e.smengoal2} holds for $\lambda > 0$ sufficiently large.
\end{lemma}

\begin{proof}
There is $I_{1/2} \in (0,\infty)$ such that the representation \eqref{e.smengoal2} of the set $[0, E_0] \cap \Sigma$ as a countable union of compact sets of small box-counting dimension is available as soon as
\begin{equation}\label{e.boxdimlemmagoal2}
\inf_{E \in \Sigma \cap [0,E_0]} I(E) > I_{1/2}.
\end{equation}

Indeed, analyticity arguments based on \cite{BLS} imply that the points of tangency between the curve of initial conditions and the center-stable manifolds of the non-wandering set of the trace map must be isolated, hence there are only finitely many those points in $\Sigma \cap [0,E_0]$, if any. If there are no points of tangency, than at every point of $\Sigma \cap [0,E_0]$ one has
$$
\dim_\mathrm{B, loc}^+(\Sigma; E)=\dim_\mathrm{H, loc}(\Sigma; E).
$$
This can be established using the same arguments that were used in the study of the spectrum of the 1D Fibonacci quantum Ising Model in \cite{Y}. Concretely, the proof of \cite[Theorem~2.1]{Y} does not use the specific features of the model, and works for the upper box counting dimension as well, as soon as the curve of initial conditions is transversal to the center-stable manifolds of the non-wandering set of the Fibonacci trace map.
The condition \eqref{e.boxdimlemmagoal2} can in turn be satisfied by choosing $\lambda$ sufficiently large. This was already discussed in the proof of \cite[Corollary~6.7]{DFG2014} and we refer the reader to that paper for the details.

In the case when there are points of tangency, denote by $A_0$ the finite set of energies in $\Sigma \cap [0,E_0]$ corresponding to points of tangency, and denote by $A_k$ the set $(\Sigma \cap [0,E_0])\backslash U_{\frac{1}{k}}(A_0)$, where $U_{\frac{1}{k}}(A_0)$ is the $\frac{1}{k}$-neighborhood of $A_0$. Then by construction, for each $k \ge 1$, at the points parameterized by energies from $A_k$, the curve of initial conditions and the center-stable manifolds are transversal, and hence the arguments above imply that $\dim_\mathrm{B}^+A_k < \frac12$ for large $\lambda$ uniformly in $k$. This proves \eqref{e.smengoal2}.
\end{proof}

With the preparatory work out of the way, we are now in a position to complete the proof of Theorem~\ref{t:lo}.

\begin{proof}[Proof of Theorem~\ref{t:lo}]
By Lemmas~\ref{l.smengoal1} and \ref{l.smengoal2}, we have \eqref{e.smengoal1} and \eqref{e.smengoal2}, which in turn imply \eqref{e.smencons1} and \eqref{e.smencons2}. Since $H_\omega \geq 0$, it follows that  $(-\infty,E_0] \cap (\Sigma+\Sigma)$ is closed, bounded, nonempty, and has empty interior. Since $\Sigma$ has no isolated points, $(-\infty,E_0] \cap (\Sigma+\Sigma)$ also has no isolated points, modulo a small wiggle of $E_0$.
\end{proof}

\begin{remark}\label{r.difficultgencase}
Notice that we referred the reader to the proof of \cite[Corollary~6.7]{DFG2014} for the argument establishing \eqref{e.boxdimlemmagoal2} in the locally constant case for sufficiently large coupling. Upon inspection of that argument, the reader will notice that any attempt at an extension to a more general setting will face difficulties of various kinds. First of all, the case of locally constant building blocks is such that solutions, and hence transfer matrix traces, and hence the invariant can be computed \emph{explicitly}. The explicit formulae obtained in this way are used in crucial ways in the proof of \cite[Corollary~6.7]{DFG2014}. Beyond the case of locally constant building blocks, one does in general not have explicit formulae for these quantities. But even if one attempted a generalization of the argument that is based on the qualitative features of these formulae, one would run into difficulties. On a technical level, the problem is related to the study of high barriers in quantum mechanics; see, for example, \cite{KMP1992, MSS1995, SS1994, S1997} for previous works on this problem. While the existing papers have investigated the reduction in tunneling upon raising a barrier, as well as resulting localization statements in the presence of a suitable sequence of high barriers, these results appear to be insufficient to produce the precise property we need from a single high barrier in order to derive a statement like \eqref{e.boxdimlemmagoal2}.
\end{remark}

\section{Discussion, Questions, and Outlook} \label{s.?}

We conclude with several questions that we view as natural follow-ups to this work, which nevertheless do not follow from our methods and require additional inspiration.

\begin{quest}\label{q.generalbumps}
Do the results of this paper hold for general $f_\a$ and $f_\b$ satisfying Assumption~\ref{assumption}? Namely, given $\ell_\alpha >0$ and $f_\alpha \in L^2(0,\ell_\alpha)$, $\alpha \in \{\a,\b\}$, is it true that $\Sigma+\Sigma$ contains a half-line? Writing $\Sigma_\lambda = \Sigma(\lambda f_\a,\lambda f_\b)$, is it true that $\Sigma_\lambda$ exhibits Cantor structure at low energies for large $\lambda$?
\end{quest}

Some of the difficulties one faces when attempting to answer Question~\ref{q.generalbumps} were discussed above in Remarks~\ref{r.difficultgencase0} and \ref{r.difficultgencase}.

\medskip

A \emph{Cantorval} is defined to be a nonempty set $S \subseteq \R$ which is compact, has dense interior, and uncountably many connected components, none of which is isolated. It is known that $C+C$ is a Cantorval for ``typical'' dynamically defined Cantor sets $C$ \cite{MMR}.

\begin{quest}\label{q.cantorval}
Is it the case that for some choice of $f_\a$ and $f_\b$, there is an interval $J \subseteq \R$ such that $J \cap (\Sigma + \Sigma)$ is a Cantorval?
\end{quest}

Clearly, the study of the structure of the sum of two Cantor sets is crucial to what we do in this paper. It is also important in several other areas and it is therefore a question that has received a lot of attention. There are well-known ways to establish that the sum is a Cantor set (e.g., by bounding the upper box dimension from above, as done in Section~\ref{s.lo}) or an interval (e.g., by bounding the thickness from below, as done in Sections~\ref{s.absc} and \ref{s.hi}). While it is known that the sum of two Cantor sets can be a Cantorval \cite{MMR, MM}, the known ways of establishing that (mostly based on \cite{MY}) 
do not seem to be applicable in our setting, where we have far too little control over the relevant properties of the 1D spectra and in which way they depend on the parameters of the model. Thus, an answer to Question~\ref{q.cantorval} appears to be currently completely out of reach.

\medskip

The Fibonacci sequence can be viewed as a coding of an irrational rotation by $\theta = \frac{\sqrt{5}-1}{2}$ via:
\[
u_n =  \begin{cases} \a & n\theta \ \mathrm{mod} \ 1 \in [1-\theta,1) \\ \b & \text{otherwise.} \end{cases}
\]
However, one can naturally replace $\theta$ by any irrational number. The resulting sequence $u$ is then called a \emph{Sturmian} sequence.
\begin{quest}
What changes if we consider potentials generated by general Sturmian sequences?
\end{quest}
In the case when the potentials are generated by a Sturmian sequence, it is still known that $\Sigma$ is a zero-measure (extended) Cantor set \cite{DFG2014}. Moreover, for typical $\theta$, there is no hyperbolic dynamical setup analogous to the one used here for $\theta = \frac{1}{2}(\sqrt{5}-1)$. Since this setup is absolutely crucial to our approach here, there is a serious obstacle to overcome in generalizing any of our results in this direction for typical $\theta$. That being said, in the event that $\theta$ is a quadratic irrationality, there is a similar (albeit more complicated) formalism \cite{G, Mei}, so it seems reasonable to expect that many of the techniques of this paper could be  applicable in that case as well.

Moreover, it is likely that the results do not in fact hold in full generality for arbitrary frequencies. For instance, in the discrete setting, it is known that the dimension of the spectrum is one for all values of the coupling constant when the frequency is very Liouville \cite{LQW}. This suggests that for Liouville $\theta$, we might not expect to see Cantor structure at the bottom of the spectrum, even for large coupling.

\begin{quest}\label{q.ap}
Can one construct $V : \R^d \to \R$ which are almost-periodic and for which the conclusion of Theorem~\ref{t.main} holds true?
\end{quest}

As mentioned earlier, there is prior work on almost periodic potentials in two dimensions for which the spectrum contains a half line \cite{KL13, KS18+}. Is it also possible to prove the Cantor structure of the spectrum at small energies? For non-separable potentials, we are lacking the tools, and for separable potentials, one would need to find 1D almost periodic potentials for which the spectrum is very thin at small energies; compare the passage from \eqref{e.smengoal2} to \eqref{e.smencons2}. In the quasi-periodic case (studied in \cite{KS18+}), no such 1D example is known, and indeed based on the known tools, it is entirely unclear how to produce one. In the limit-periodic case (studied in \cite{KL13}), there is a dichotomy, which is a function of the existing methods. If the potentials are such that the existence of a half line in the spectrum of the 2D model can be proved, the rate of approximation by periodic potentials is so fast that at small energies, the 1D model has spectrum that is too thick for any known technique to establish the Cantor structure of the sum set. If on the other hand, the rate of approximation is such that the 1D spectrum is sufficiently thin at small energies such that the 2D spectrum has a Cantor structure there, then the spectrum is very thin at all energies and the 2D spectrum is nowhere dense throughout. This mechanism underlies the work \cite{DFG2019}.

The dichotomy above, showing that one cannot have both desirable properties at the same time, is related to the known phenomenon that for discrete 1D Schr\"odinger operators with limit-periodic potentials, the spectral type and the qualitative features of the spectrum are the same for all positive values of the coupling constant in all examples in which one is able to prove statements for more than one value of the coupling constant (cf.~\cite{DF2019+}). This translates to the rough energy-independence of the spectral features discussed above for continuum 1D Schr\"odinger operators with limit-periodic potentials.

\medskip

While we focused on the topological structure of the spectrum in the present paper, the following question is naturally of interest as well:

\begin{quest}
What can be said about the type of the spectral measures and the density of states measure, both for the specific operators studied in this paper and more generally for multi-dimensional Schr\"odinger operators with almost periodic or uniformly recurrent potential.
\end{quest}

The papers \cite{KL13, KS18+} mentioned before also address the spectral type in the high-energy region and prove that it is purely absolutely continuous under their assumptions. Those papers are unable to say anything about the spectral type in the low or medium-energy regime. On the other hand, the operators studied in \cite{DFG2019} have purely singular continuous spectrum: the absence of eigenvalues follows from the Gordon lemma, applied to the underlying 1D operators, and the preservation of that fact when forming separable potentials from these 1D pieces; while the absence of absolutely continuous spectrum is immediate from the fact that the spectrum has zero Lebesgue measure. For the operators studied in the present paper, the singular continuous nature of the spectral measure and the density of states measure in the low-energy region follows again from the same arguments, namely Gordon lemma and zero-measure spectrum.\footnote{While \cite{DFG2014} did not make the application of the Gordon lemma to the 1D operators studied there explicit, the argument is completely analogous to the one used in the discrete case \cite{DL99} to exclude eigenvalues.}
We are unable to determine the spectral type of these operators in the large-energy region. Based on analogy with the work \cite{DGS} on the weakly-coupled discrete case, we would expect that at least the density of states measure is purely absolutely continuous there, but more work will be necessary to prove this. Notice that for similar models, the presence of an absolutely continuous component in the density of states measure in the large-energy region was established in \cite{FTY}. It is also worth mentioning that the relation between the structure of the spectrum as a set and properties of the density of states measure can be quite non-trivial. For example, in the case of the discrete square Fibonacci Hamiltonian there are regimes with positive measure spectrum but singular density of states measures, see \cite{DG2018}. Similar mechanisms could be present in the continuum case as well. Whether the spectral measures are purely absolutely continuous in the high-energy region for the operators studied here is currently wide open; proving this seems to be well beyond current methods.

Finally we want to mention that, while they do not determine the type of the measures in question, Parnovski and Shterenberg study them in the high-energy regime for quite general multi-dimensional almost-periodic Schr\"odinger operators in the papers \cite{PS2012, PS2016}, where they obtain complete asymptotic expansions for both the integrated density of states and the spectral function.

\end{document}